\PassOptionsToPackage{unicode}{hyperref}
\PassOptionsToPackage{hyphens}{url}
\PassOptionsToPackage{dvipsnames,svgnames,x11names}{xcolor}
\documentclass[
  12pt]{article}

\usepackage{amsthm}
\usepackage{bbm}
\usepackage{epsfig}
\usepackage{graphicx}
\usepackage{mathrsfs}
\usepackage{amsfonts}
\usepackage{color}
\usepackage{rotating}
\usepackage{verbatim}
\usepackage{epstopdf}
\usepackage[authoryear,round,comma]{natbib}
\usepackage{hyperref}
\usepackage{mathtools}
\usepackage{amssymb}
\usepackage{pifont}

\newtheorem{lemma}{Lemma}[section]

\newtheorem{theorem}{Theorem}[section]
\newtheorem{corollary}{Corollary}[section]

\usepackage{amsmath,amssymb}
\usepackage{iftex}
\ifPDFTeX
  \usepackage[T1]{fontenc}
  \usepackage[utf8]{inputenc}
  \usepackage{textcomp} 
\else 
  \usepackage{unicode-math}
  \defaultfontfeatures{Scale=MatchLowercase}
  \defaultfontfeatures[\rmfamily]{Ligatures=TeX,Scale=1}
\fi
\usepackage{lmodern}
\ifPDFTeX\else
\fi
\IfFileExists{upquote.sty}{\usepackage{upquote}}{}
\IfFileExists{microtype.sty}{
  \usepackage[]{microtype}
  \UseMicrotypeSet[protrusion]{basicmath} 
}{}
\makeatletter
\@ifundefined{KOMAClassName}{
  \IfFileExists{parskip.sty}{%
    \usepackage{parskip}
  }{
    \setlength{\parindent}{0pt}
    \setlength{\parskip}{6pt plus 2pt minus 1pt}}
}{
  \KOMAoptions{parskip=half}}
\makeatother
\usepackage{xcolor}
\makeatletter
\ifx\paragraph\undefined\else
  \let\oldparagraph\paragraph
  \renewcommand{\paragraph}{
    \@ifstar
      \xxxParagraphStar
      \xxxParagraphNoStar
  }
  \newcommand{\xxxParagraphStar}[1]{\oldparagraph*{#1}\mbox{}}
  \newcommand{\xxxParagraphNoStar}[1]{\oldparagraph{#1}\mbox{}}
\fi
\ifx\subparagraph\undefined\else
  \let\oldsubparagraph\subparagraph
  \renewcommand{\subparagraph}{
    \@ifstar
      \xxxSubParagraphStar
      \xxxSubParagraphNoStar
  }
  \newcommand{\xxxSubParagraphStar}[1]{\oldsubparagraph*{#1}\mbox{}}
  \newcommand{\xxxSubParagraphNoStar}[1]{\oldsubparagraph{#1}\mbox{}}
\fi
\makeatother

\usepackage{longtable,booktabs,array}
\usepackage{calc} 
\usepackage{etoolbox}
\makeatletter
\patchcmd\longtable{\par}{\if@noskipsec\mbox{}\fi\par}{}{}
\makeatother
\IfFileExists{footnotehyper.sty}{\usepackage{footnotehyper}}{\usepackage{footnote}}
\makesavenoteenv{longtable}
\usepackage{graphicx}
\makeatletter
\def\maxwidth{\ifdim\Gin@nat@width>\linewidth\linewidth\else\Gin@nat@width\fi}
\def\maxheight{\ifdim\Gin@nat@height>\textheight\textheight\else\Gin@nat@height\fi}
\makeatother
\setkeys{Gin}{width=\maxwidth,height=\maxheight,keepaspectratio}
\makeatletter
\def\fps@figure{htbp}
\makeatother

\makeatletter
\@ifpackageloaded{caption}{}{\usepackage{caption}}
\AtBeginDocument{%
\ifdefined\contentsname
  \renewcommand*\contentsname{Table of contents}
\else
  \newcommand\contentsname{Table of contents}
\fi
\ifdefined\listfigurename
  \renewcommand*\listfigurename{List of Figures}
\else
  \newcommand\listfigurename{List of Figures}
\fi
\ifdefined\listtablename
  \renewcommand*\listtablename{List of Tables}
\else
  \newcommand\listtablename{List of Tables}
\fi
\ifdefined\figurename
  \renewcommand*\figurename{Figure}
\else
  \newcommand\figurename{Figure}
\fi
\ifdefined\tablename
  \renewcommand*\tablename{Table}
\else
  \newcommand\tablename{Table}
\fi
}
\@ifpackageloaded{float}{}{\usepackage{float}}
\floatstyle{ruled}
\@ifundefined{c@chapter}{\newfloat{codelisting}{h}{lop}}{\newfloat{codelisting}{h}{lop}[chapter]}
\floatname{codelisting}{Listing}

\makeatother
\makeatletter
\@ifpackageloaded{caption}{}{\usepackage{caption}}
\@ifpackageloaded{subcaption}{}{\usepackage{subcaption}}
\makeatother

\ifLuaTeX
  \usepackage{selnolig}  
\fi
\usepackage[]{natbib}
\usepackage{bookmark}

\IfFileExists{xurl.sty}{\usepackage{xurl}}{} 
\hypersetup{
  pdftitle={Title},
  pdfauthor={Author 1; Author 2},
  pdfkeywords={3 to 6 keywords, that do not appear in the title},
  colorlinks=true,
  linkcolor={blue},
  filecolor={Maroon},
  citecolor={Blue},
  urlcolor={Blue},
  pdfcreator={LaTeX via pandoc}}

\newcommand{\anon}{1}

\begin{document}

\def\spacingset#1{\renewcommand{\baselinestretch}%
{#1}\small\normalsize} \spacingset{1}


\if1\anon
{
  \title{\bf A Generative High Quantile Homogeneity Test Using Bahadur Representation for Heteroskedastic High Quantile Regression of Tail Dependent Time Series}
  \author{Ting Zhang\thanks{
    Ting Zhang is supported by the NSF Grant DMS-2412661.}\hspace{.2cm}\\
    Department of Statistics, University of Georgia\\
    and \\
    Fangwei Wu\\
    Department of Statistics and Data Science, Tsinghua University\\
    and \\
    Jingying Gao\\
    Department of Statistics, University of Georgia}
  \maketitle
} \fi

\if0\anon
{
  \bigskip
  \bigskip
  \bigskip
  \begin{center}
    {\LARGE\bf A Generative High Quantile Homogeneity Test Using Bahadur Representation for Heteroskedastic High Quantile Regression of Tail Dependent Time Series}
\end{center}
  \medskip
} \fi

\bigskip
\begin{abstract}
We consider a high quantile homogeneity test to determine whether a certain set of explanatory variables has homogeneous effects on different high quantiles of the response variable in the tail. To accommodate for situations under both the null and the alternative, the auxiliary process in this case may no longer be treated as stationary, and the problem requires a joint analysis of both homoscedastic and heteroskedastic high quantiles. For this, we develop a novel Bahadur representation result in the high quantile setting for a general class of tail dependent time series under potential heteroskedasticity, which can be of interest by its own. In particular, the Bahadur representation provides a foundation for reducing problems regarding nonlinear high quantile regression estimators to those regarding suitably constructed linear forms with an explicit error bound and can be transformative and useful in many statistical problems. We in the current article apply it to guide the development of a generative high quantile homogeneity test, which is then illustrated through applications to both synthetic and real data.
\end{abstract}

\noindent%
\vfill

\newpage
\spacingset{1.8} 

\section{Introduction}\label{sec:introduction}
In many applications it is desirable to study and compare the effect of a certain set of regressors on different tail quantiles of the response. For example, \citet{Elsner:Kossin:Jagger:2008} fitted linear trends for upper high quantiles of a satellite-derived data that contains lifetime-maximum wind speeds of tropical cyclones, and they compared the linear coefficients obtained from quantile regression models at different quantile levels to study the trend difference of the strongest tropical cyclones. \citet{Zhang:2021} considered fitting polynomial trends for upper high quantiles of the global temperature anomaly series, and it was identified that different polynomials should be used for trend modeling of the 95\% and 99\% quantiles, indicating a different temperature pattern of the hottest days when compared with the average. The comparisons in the aforementioned analyses, however, were mainly conducted heuristically using standard errors calculated from separate quantile regression models assuming either independent data or stationary homoscedastic auxiliary processes. To provide a rigorous homogeneity test that can survive under the dependence between different high quantiles, however, it requires a joint analysis of both homoscedastic and heteroskedastic high quantiles under the null and the alternative, and as a result the auxiliary process can no longer be treated as stationary. For this, we develop a novel Bahadur representation in the high quantile setting for a general class of tail dependent time series so that the difficult problem of comparing nonlinear high quantile regression estimators with an unknown dependence structure can be translated into a simpler problem about comparing suitably constructed linear tail quantities.

In the context of quantile regression \citep{Koenker:2005}, the Bahadur representation refers to the use of a linear form to represent the originally nonlinear quantile regression estimator with an explicit approximation bound on the difference between the two so that the result can be generally useful in different contexts. The problem of establishing the Bahadur representation, however, can often be quite involved and has attracted a considerable attention in the literature; see for example \citet{Bahadur:1966}, \citet{Kiefer:1967}, \citet{Hesse:1990}, \citet{He:Shao:1996}, \citet{Ho:Hsing:1996}, \citet{Wu:2005:Bahadur}, \citet{Wu:2007:MEstimation}, \citet{Wu:Zhou:2017}, \citet{Wu:Yu:Wang:2021}, \citet{Berrahou:Bouzebda:Douge:2024}, and references therein. Nevertheless, existing results in this direction were mostly developed in the conventional quantile regression setting where the quantile level is treated as fixed, and as a result they may not be directly applicable or suitable for tail analyses that focus more on the study of rare or extremal events. We in this article fill the gap by developing a general Bahadur representation result for tail dependent data with possible heteroskedasticity under the double asymptotic scheme where the quantile level is no longer treated as fixed but allowed to approach the boundary as the sample size increases. Similar to the role of an increasing dimension in high-dimensional problems, the double asymptotic scheme in the current high quantile setting may also alter the asymptotic behavior and makes it more challenging to study the Bahadur representation in high quantile regression models. It can be seen from Section \ref{sec:Bahadur} that, compared with its fixed quantile counterpart, the Bahadur representation in the high quantile setting involves a linear approximation with different rates, which can depend on an interplay between the degree of tail extremeness and the tail behavior of the underlying distribution. In addition, unlike \citet{Zhang:2021} which could afford to perform the martingale decomposition directly on the linear form, the current Bahadur representation problem requires a martingale decomposition on the original nonlinear function in order to quantify the explicit distance between the nonlinear high quantile regression estimator and its linear approximation. To handle this, we first use the inequality of \citet{Freedman:1975} to quantify the pointwise tail probability of martingales and then combine it with a deep probabilistic chaining argument to obtain a uniform result on an uncountable set that was otherwise not needed in \citet{Zhang:2021}. It is remarkable that the Bahadur representation may be of interest by its own, as it provides a foundational tool for reducing problems regarding nonlinear high quantile regression estimators to those regarding suitably constructed linear forms with an explicit error bound and can be transformative and useful in many statistical problems.

We in the current article apply the newly developed Bahadur representation to guide the development of a high quantile homogeneity test. \citet{Wang:Li:He:2012} considered the problem of extrapolating high quantiles using quantile regression, and they proposed two different estimation methods depending on whether the regressor effect is believed to be homogeneous or not. Since the two methods generally lead to different estimates, to address the ambiguity, \citet{Wang:Li:He:2012} recommended using the analysis of variance (ANOVA) test of \citet{Koenker:Bassett:1982} for heteroskedasticity to decide which estimator to use. Nevertheless, the test of \citet{Koenker:Bassett:1982} was developed for independent data under the conventional quantile regression setting when the quantile levels are assumed to be fixed. In time series applications, however, dependence is the rule rather than the exception and it has been demonstrated in the literature that ignoring the dependence can lead to erroneous conclusions in quantile analyses; see for example \citet{Wu:2007:MEstimation}, \citet{Zhou:Shao:2013} and \citet{Zhang:2021}. It is remarkable that dealing with dependent data in high quantile regression problems can be a quite nontrivial task. For example, \citet{Chernozhukov:2005} considered the use of the strong mixing condition \citep{Rosenblatt:1956} while still having to impose an additional condition to limit the dependence in the tail to a negligible degree for asymptotic theory of tail quantities; see also \citet{Chernozhukov:FernandezVal:2011} for the use of a similar negligibility condition. \citet{Daouia:Stupfler:Usseglio-Carleve:2023} considered the situation when data points used in each particular estimation are well separated by a state variable and can thus be treated as asymptotically independent. We in this article provide a rigorous statistical test on the homogeneity of regressor effect on different high quantiles for a general class of tail dependent time series, and it can be seen from the derived limiting distribution that the asymptotic behavior of the test statistic can indeed be affected by the existence of dependence in the tail region. The test aims to determine statistically if a certain set of explanatory variables have homogeneous effects on different high quantiles of the response variable in the tail, and we consider both situations when the quantile levels to be compared approach the extreme at different rates and when they share the same rate of extremeness but with different constants. It can be seen from our developed asymptotic theory that the two situations are guided by different limit theorems for inference, while in practice it can be ambiguous to determine which limit theorem to use in a given application. For this, we propose a generative high quantile regression approach that is able to not only unify the two situations without having to deal with the annoying ambiguity but also lead to improved finite-sample performance of the developed test.

The rest of the article is organized as follows. Section \ref{sec:framework} introduces the mathematical setting and framework. The high quantile Bahadur representation result is developed in Section \ref{sec:Bahadur} for a general class of tail dependent time series, which is then applied in Section \ref{sec:application} to guide the development of a generative high quantile homogeneity test. Numerical experiments including Monte Carlo simulations and a real data application are also provided to illustrate the proposed test. Section \ref{sec:discussion} concludes the article by providing a discussion.

\section{Setting and Framework}\label{sec:framework}
Consider the quantile regression model \citep{Koenker:Bassett:1978}
\begin{equation}\label{eqn:highquantileregressionmodel}
	y_{i,n} = \boldsymbol x_{i,n}^\top \boldsymbol \beta_{1-\alpha_n,n} + \tilde e_{i,n},\quad i = 1,\ldots,n,
\end{equation}
where $\boldsymbol x_{i,n}^\top \boldsymbol \beta_{1-\alpha_n,n}$ represents the $(1-\alpha_n)$-th quantile with $^\top$ denoting the matrix transpose for some regression coefficient $\boldsymbol \beta_{1-\alpha_n,n} \in \mathbb R^p$ and $\tilde e_{i,n} = y_{i,n} - \boldsymbol x_{i,n}^\top \boldsymbol \beta_{1-\alpha_n,n}$ is the associated auxiliary variable satisfying $\mathrm{pr}(\tilde e_{i,n} \leq 0) = \mathrm{pr}(y_{i,n} \leq \boldsymbol x_{i,n}^\top \boldsymbol \beta_{1-\alpha_n,n}) = 1-\alpha_n$. In the current high quantile regression setting, the quantile level $1-\alpha_n \to 1$ as $n \to \infty$ is under the double asymptotic scheme and the auxiliary variables $(\tilde e_{i,n})$ form a triangular array. Intuitively, (\ref{eqn:highquantileregressionmodel}) decomposes $y_{i,n}$ into its $(1-\alpha_n)$-th quantile $\boldsymbol x_{i,n}^\top \boldsymbol \beta_{1-\alpha_n,n}$ and a remainder $\tilde e_{i,n}$; see also the discussions in \citet{Bai:Rao:Wu:1992}, \citet{Wu:2007:MEstimation} and \citet{Zhang:2021}. For fixed quantile levels, \citet{Bai:Rao:Wu:1992} modeled the auxiliary variables $\tilde e_{i,n}$, $i = 1,\ldots,n$, as independent and identically distributed (i.i.d.), while \citet{Wu:2007:MEstimation} extended their results by allowing $(\tilde e_{i,n})$ to form a stationary time series. Recently, \citet{Zhang:2021} explored the high quantile regression setting when $1-\alpha_n \to 1$ as $n \to \infty$ by assuming that the associated high quantile auxiliary variables $\tilde e_{i,n}$, $i = 1,\ldots,n$, are strictly stationary. When the regressor effect is not homogeneous as in the non-common-slope model of \citet{Wang:Li:He:2012}, however, the auxiliary variables $\tilde e_{i,n}$, $i = 1,\ldots,n$, can no longer be treated as strictly stationary as the inhomogeneous regressor effect can cause heteroskedasticity in the auxiliary variables. This motivates us to consider a more general class of tail dependent auxiliary processes
\begin{equation}\label{eqn:eintildedefinition}
	\tilde e_{i,n} = \gamma_{i,n} e_{i,n}
\end{equation}
for some $0 < \gamma_{i,n} < \infty$ that can possibly depend on $\boldsymbol x_{i,n}$ and
\begin{equation*}
	e_{i,n} = G_n(\mathcal F_i),\quad \mathcal F_i = (\ldots,\epsilon_{i-1},\epsilon_i),
\end{equation*}
where $\epsilon_k$, $k \in \mathbb Z$, are i.i.d. innovations and $G_n$ is a sequence of measurable functions such that $e_{i,n}$ is properly defined. Following \citet{Wu:2005}, we interpret $e_{i,n}$ as the output of the physical system $G_n$ with input $\mathcal F_i$ in a triangular array setting, which is quite general and covers a huge class of popular time series models as special cases; see also the discussions in \citet{Wiener:1958}, \citet{Tong:1990}, \citet{Zhang:2005}, \citet{Wu:2011}, and \citet{Zhang:Xu:2024}. We shall here in addition consider a generalization by allowing the auxiliary process to exhibit heteroskedasticity as in (\ref{eqn:eintildedefinition}). Since the decomposition of $\gamma_{i,n} e_{i,n}$ in (\ref{eqn:eintildedefinition}) is not unique, without loss of generality we assume throughout the article that $n^{-1} \sum_{i=1}^n \gamma_{i,n} = 1$ to ensure the identifiability.

Let $\epsilon_0^\star$ be identically distributed as $\epsilon_0$ but independent of $(\epsilon_k)_{k \in \mathbb Z}$, then we can define the coupled filtration $\mathcal F_k^\star = (\mathcal F_{-1},\epsilon_0^\star,\epsilon_1,\ldots,\epsilon_k)$ and its associated output $e_{k,n}^\star = G_n(\mathcal F_k^\star)$. Therefore, $e_{k,n}^\star$ represents the output of the same physical system but with the innovation at time zero $\epsilon_0$ replaced by its independent copy $\epsilon_0^\star$ in the input filtration. Let $F_n(u) = \mathrm{pr}(e_{k,n} \leq u)$ denote the marginal distribution function of $e_{k,n}$, then we follow \citet{Zhang:2021} and consider
\begin{equation*}
	\theta_{n,\alpha}(k) = \sup_{a \in (0,\alpha], N \geq n} \mathrm{pr}\{e_{k,N}^\star \leq F_N^{-1}(1-a) \mid e_{k,N} > F_N^{-1}(1-a)\},\quad \alpha \in (0,1),
\end{equation*}
which measures the degree of tail dependence through the effect of tail adversarial coupling; see also the discussions in \citet{Zhang:2022} and \citet{Bai:Zhang:2024}. Write
\begin{equation*}
	\Theta_{n,\alpha,q}(m) = \sum_{k=m}^\infty \{\theta_{n,\alpha}(k)\}^{1/q},\quad m \geq 0,
\end{equation*}
which measures the cumulative tail adversarial effect of the current innovation on all future observations, then the process $(e_{i,n})$ is said to be asymptotically tail adversarial $q$-stable, or $(e_{i,n}) \in \mathrm{TAS}_q$, if
\begin{equation*}
	\lim_{\alpha \downarrow 0} \lim_{n \to \infty} \Theta_{n,\alpha,q}(0) < \infty.
\end{equation*}
As illustrated in \citet{Zhang:2021}, \citet{Zhang:2022}, \citet{Bai:Zhang:2024}, \citet{Zhang:Xu:2024} and \citet{Cao:Gao:Shao:Sriram:Wang:Wen:Zhang:2025}, the tail adversarial stability framework can often lead to weaker conditions than the conventional strong mixing framework \citep{Rosenblatt:1956,Chernozhukov:2005} and be useful in obtaining unprecedented limit theorems of nontrivial tail quantities that could not have been thoroughly explored otherwise. We shall here provide a further illustration of the aforementioned tail adversarial stability framework using the moving-maximum process of \citet{Hall:Peng:Yao:2002}, which was shown to be dense in the class of stationary processes whose finite-dimensional distributions are extreme-value of a given type; see also the discussion in \citet{Zhang:Zhang:Cui:2017}. For this, let $\epsilon_k$, $k \in \mathbb Z$, be independent Fr\'{e}chet random variables with distribution function $F_\epsilon(z) = \mathrm{pr}(\epsilon_k \leq z) = \exp(-z^{-\kappa})$ for some $\kappa > 0$, and we consider the moving-maximum process
\begin{equation}\label{eqn:einmm}
	e_{i,n} = \max_{0 \leq l < \infty} a_l \epsilon_{i-l} - Q_{e,1-\alpha_n},
\end{equation}
which is well-defined if the nonnegative coefficients satisfy $\sum_{l=0}^\infty a_l^\kappa < \infty$ and $Q_{e,1-\alpha_n}$ is the $(1-\alpha_n)$-th quantile of $\max_{0 \leq l < \infty} a_l \epsilon_{i-l}$ so that $\mathrm{pr}(\gamma_{i,n} e_{i,n} \leq 0) = \mathrm{pr}(e_{i,n} \leq 0) = 1-\alpha_n$ represents a proper quantile regression auxiliary process; see for example \citet{Hall:Peng:Yao:2002}, \citet{Zhang:Smith:2004}, \citet{Heffernan:Tawn:Zhang:2007}, \citet{Zhang:Zhang:Cui:2017} and \citet{Zhang:2021}. For the moving-maximum process (\ref{eqn:einmm}), by elementary calculations one can show that
\begin{equation*}
	\theta_{n,\alpha}(k) \leq 2 \left(\sum_{l=0}^\infty a_l^\kappa\right)^{-1} a_k^\kappa
\end{equation*}
holds for any $\alpha \in (0,1/2)$, and as a result the process $(e_{i,n}) \in \mathrm{TAS}_q$ if $\sum_{l=0}^\infty a_l^{\kappa/q} < \infty$. Therefore, in view of the existence condition $\sum_{l=0}^\infty a_l^\kappa < \infty$ under which the process is well-defined, the tail adversarial stability condition seems to be mild and reasonable. We shall here take advantage of the tail adversarial stability framework to develop a novel Bahadur representation result with an explicit approximation bound for high quantile regression estimators in the challenging double asymptotic scheme when the quantile level is allowed to approach the extreme as the sample size grows to infinity for a general class of processes that can exhibit dependence in both tail and non-tail regions with possible heteroskedasticity.

\section{High Quantile Bahadur Representation}\label{sec:Bahadur}
To estimate the quantile regression model (\ref{eqn:highquantileregressionmodel}) when the quantile level $1-\alpha_n \to 1$ as $n \to \infty$, we consider the high quantile regression estimator
\begin{equation*}
	\hat{\boldsymbol \beta}_{1-\alpha_n,n} = \mathop{\mathrm{argmin}}_{\boldsymbol \beta \in \mathbb R^p} \sum_{i=1}^n \rho_{1-\alpha_n}(y_{i,n} - \boldsymbol x_{i,n}^\top \boldsymbol \beta),
\end{equation*}
where $\rho_{1-\alpha_n}(u)=(1-\alpha_n)u^+ + \alpha_n(-u)^+$ is the check function with $u^+ = \max(u,0)$ and we denote its left derivative by $\psi_{1-\alpha_n}(u) = (1-\alpha_n) - \mathbbm{1}_{\{u \leq 0\}}$ with $\mathbbm{1}_{\{\cdot\}}$ being the indicator function. Assuming that $\Sigma_n = n^{-1} \sum_{i=1}^n \boldsymbol x_{i,n} \boldsymbol x_{i,n}^\top$ is non-singular for all large $n$, and we define $\boldsymbol z_{i,n} = \Sigma_n^{-1/2} \boldsymbol x_{i,n}$. Then $n^{-1}\sum_{i=1}^n \boldsymbol z_{i,n} \boldsymbol z_{i,n}^\top = \mathrm{I}_{p \times p}$, the $p \times p$ identity matrix, and it becomes more convenient to consider the rescaled model
\begin{equation*}
	y_{i,n} = \boldsymbol z_{i,n}^\top \boldsymbol \varphi_{1-\alpha_n,n} + \gamma_{i,n} e_{i,n},\quad i = 1,\ldots,n,
\end{equation*}
which is mathematically equivalent to the heteroskedastic high quantile regression model (\ref{eqn:highquantileregressionmodel}) with $\boldsymbol \varphi_{1-\alpha_n,n} = \Sigma_n^{1/2} \boldsymbol \beta_{1-\alpha_n,n}$. The associated high quantile regression estimator is then given by
\begin{equation*}
	\hat{\boldsymbol \varphi}_{1-\alpha_n,n} = \mathop{\mathrm{argmin}}_{\boldsymbol \varphi \in \mathbb R^p} \sum_{i=1}^n \rho_{1-\alpha_n}(y_{i,n} - \boldsymbol z_{i,n}^\top \boldsymbol \varphi),
\end{equation*}
which satisfies the relation $\hat{\boldsymbol \varphi}_{1-\alpha_n,n} = \Sigma_n^{1/2} \hat{\boldsymbol \beta}_{1-\alpha_n,n}$. The major goal of this section is to develop a general Bahadur representation result for the high quantile regression estimator so that we can write the nonlinear estimator $\hat{\boldsymbol \varphi}_{1-\alpha_n,n}$ as the sum of a suitably constructed linear form and a remainder term with an explicit order. For this, we make the following assumptions.
\begin{itemize}
	\item [(A1)] The triangular array $(e_{i,n})\in \mathrm{TAS}_q$ for some $q \geq 2$.
	\item [(A2)] There exists an $\alpha\in(0,1)$ such that $F_n(\cdot)$ is continuously differentiable with uniformly bounded and strictly positive derivative $f_n(\cdot)$ in its upper tail $\{F_n^{-1}(1-\alpha),F_n^{-1}(1)\}$ with $
	\lim\inf_{n\rightarrow \infty}|F_n^{-1}(1)-F_n^{-1}(1-\alpha)|>0$ for all large $n$.
	\item [(A3)] The rescaled design satisfies $\max_{1\leq i\leq n} |\boldsymbol z_{i,n}| = o\{(n\alpha_n)^{1/2}\}$.
	\item [(A4)] There exists $0 < c_{\gamma,1} \leq c_{\gamma,2} < \infty$ such that
	\begin{equation*}
		c_{\gamma,1} \leq \min_{1 \leq i \leq n} \gamma_{i,n} \leq \max_{1 \leq i \leq n} \gamma_{i,n} \leq c_{\gamma,2}
	\end{equation*}
	holds for all large $n$.
\end{itemize}

We shall here provide a brief discussion on the aforementioned conditions. In particular, condition (A1) refers to the tail adversarial stability discussed in Section \ref{sec:framework}, which has been studied for various tail dependent processes in \citet{Zhang:2021}, \citet{Zhang:2022}, \citet{Bai:Zhang:2024} and \citet{Zhang:Xu:2024}. Condition (A2) concerns the tail smoothness of the underlying distribution, which is mild and satisfied by many common distributions such as the uniform, normal, exponential and heavy-tailed Pareto distributions. Condition (A3) is a mild Lindeberg-type condition on the rescaled design in the current high quantile regression setting, which has been used and verified in \citet{Zhang:2021}. The aforementioned article studied the convergence rate and central limit theorem for high quantile regression estimators when the auxiliary process is strictly stationary. We shall here consider the more challenging Bahadur representation problem under the more general setting when the auxiliary process is not necessarily stationary with potential heteroskedasticity. Condition (A4) is natural in the context of heteroskedastic quantile regression models; see for example \citet{Wang:Li:He:2012}. To establish the Bahadur representation, we introduce the triangular array function
\begin{equation*}
	\boldsymbol K_n(\boldsymbol \varphi) = \boldsymbol \Omega_n(\boldsymbol \varphi) - E\{\boldsymbol \Omega_n(\boldsymbol \varphi)\},\quad \boldsymbol \Omega_n(\boldsymbol \varphi) = \sum_{i=1}^n \psi_{1-\alpha_n}(\gamma_{i,n} e_{i,n} - \boldsymbol z_{i,n}^\top \boldsymbol \varphi)\boldsymbol z_{i,n},
\end{equation*}
which relates to the derivative of the high quantile regression criterion function. Let
\begin{equation*}
	\boldsymbol M_n(\boldsymbol \varphi) =  \sum_{i=1}^n\left[\psi_{1-\alpha_n}(\gamma_{i,n} e_{i,n} - \boldsymbol z_{i,n}^\top \boldsymbol \varphi) -E\{\psi_{1-\alpha_n}(\gamma_{i,n} e_{i,n} - \boldsymbol z_{i,n}^\top \boldsymbol \varphi) \mid \mathcal{F}_{i-1} \}\right]\boldsymbol z_{i,n}
\end{equation*}
be the martingale generated from the nonlinear summands in $\boldsymbol \Omega_n(\boldsymbol \varphi)$, then we can write
\begin{equation*}
	\boldsymbol K_n(\boldsymbol \varphi) =  \boldsymbol M_n(\boldsymbol \varphi) + \boldsymbol N_n(\boldsymbol \varphi),
\end{equation*}
where
\begin{equation*}
	\boldsymbol N_n(\boldsymbol \varphi) =  \sum_{i=1}^n\left[E\{\psi_{1-\alpha_n}(\gamma_{i,n} e_{i,n} - \boldsymbol z_{i,n}^\top \boldsymbol \varphi) \mid \mathcal{F}_{i-1} \} - E\{\psi_{1-\alpha_n}(\gamma_{i,n} e_{i,n} - \boldsymbol z_{i,n}^\top \boldsymbol \varphi)\}\right] \boldsymbol z_{i,n}
\end{equation*}
is used to denote the remainder term. Unlike \citet{Zhang:2021} which had the ease of being able to perform the martingale decomposition directly on a linear form, the current Bahadur representation requires a martingale decomposition on the original nonlinear function in order to explicitly characterize the difference between the nonlinear high quantile regression estimator and its suitably constructed linear approximation. As a result, the theoretical derivation is much more involved and uses a different set of techniques. In particular, we use the martingale inequality of \citet{Freedman:1975} to quantify the pointwise tail probability and then combine it with a deep probabilistic chaining argument to obtain a uniform result on an uncountable set that was otherwise not needed in \citet{Zhang:2021}. Let $\boldsymbol 0_p$ denote the $p$-dimensional vector with zero components and write $r_n = \max_{1 \leq i \leq n}|\boldsymbol z_{i,n}|$. Lemma \ref{lem:Mn} provides a uniform bound for the martingale part in the high quantile setting, which is useful in obtaining the Bahadur representation.

\begin{lemma}\label{lem:Mn}
	Assume (A4), $\alpha_n\to 0$ and $n\alpha_n\to  \infty$. If $\iota_{1,n}(u; \mathcal{F}_0) = \mathrm{pr}(e_{1,n} > u \mid \mathcal F_0)$ is stochastically Lipschitz continuous in a neighborhood of zero in the sense that there exists an $\varepsilon_0 > 0$ such that
	\begin{equation}\label{eqn:L0ncondition}
		L_{0,n} = \sup_{-\varepsilon_0 \leq u_1 < u_2 \leq \varepsilon_0}{|\iota_{1,n}(u_2; \mathcal F_0) - \iota_{1,n}(u_1; \mathcal F_0)| \over |u_2-u_1|} \in \mathcal L^1,
	\end{equation}
	then for any sequence $\delta_n > 0$ with $\delta_n r_n = o(1)$ we have
	\begin{equation*}
		\sup_{|\boldsymbol \varphi| \leq \delta_n}|\boldsymbol M_n(\boldsymbol \varphi) - \boldsymbol M_n(\boldsymbol 0_p)| = O_p[\{\varrho_n(\delta_n)\}^{1/2}\log n + r_n^2 n^{-3}],
	\end{equation*}
	where $\varrho_n(\delta)=\sum_{i=1}^n|\boldsymbol z_{i,n}|^2\{F_n(c_{\gamma,1}^{-1}|\boldsymbol z_{i,n}|\delta) - F_n(-c_{\gamma,1}^{-1}|\boldsymbol z_{i,n}|\delta)\}$.
\end{lemma}

Condition (\ref{eqn:L0ncondition}) concerns the smoothness of the tail conditional distribution and is in general very mild. In particular, if we use $F_{1,n}(u \mid \mathcal F_0) = \mathrm{pr}(e_{1,n} \leq u \mid \mathcal F_0)$ to denote the conditional distribution function, then condition (\ref{eqn:L0ncondition}) is satisfied if the associated conditional density $f_{1,n}(u \mid \mathcal F_0) = \partial F_{1,n}(u \mid \mathcal F_0)/\partial u$ is bounded in a neighborhood of zero. Note that zero serves as the $(1-\alpha_n)$-th quantile of the auxiliary process and thus falls into the tail region. Similarly, write $\iota_{1,n}(u; \mathcal F_i) = \mathrm{pr}(e_{i+1,n} > u \mid \mathcal F_i)$ and $\iota_{1,n}^{(q)}(u; \mathcal F_{i}) = \partial^q \iota_{1,n}(u; \mathcal F_{i})/\partial u^q$, then Lemma \ref{lem:Nn} provides a bound for the remainder term.

\begin{lemma}\label{lem:Nn}
	Assume (A4), $\alpha_n\to 0$ and $n\alpha_n\to  \infty$. If there exists an $\varepsilon_0 > 0$ such that
    \begin{equation*}
    \sup_{|\varepsilon| \leq \varepsilon_0} \|\iota_{1,n}^{(q)}(\varepsilon; \mathcal F_{i})\| < \infty
    \end{equation*}
    and
	\begin{equation}\label{eqn:iotacondition}
		\sum_{i=0}^\infty \sup_{|\varepsilon| \leq \varepsilon_0}\|\mathrm{E}\{\iota_{1,n}^{(q)}(\varepsilon; \mathcal{F}_{i}) \mid \mathcal{F}_{0}\} - \mathrm{E}\{\iota_{1,n}^{(q)}(\varepsilon; \mathcal{F}_{i}^\star) \mid \mathcal{F}_{0}^\star\}\| = O(\alpha_n^{1/2})
	\end{equation}
	hold for $q = 1,\ldots,p$, then for any sequence $\delta_n > 0$ with $\delta_n r_n = o(1)$ we have
	\begin{equation*}
		\left\|\sup_{|\boldsymbol \varphi| \leq \delta_n}|\boldsymbol N_n(\boldsymbol \varphi) - \boldsymbol N_n(\boldsymbol 0_p)|\right\| = O_p[\delta_n \{\alpha_n \nu_n(4)\}^{1/2}],
	\end{equation*}
	where $\nu_n(q) = \sum_{i=1}^n|\boldsymbol z_{i,n}|^q$.
\end{lemma}

Condition (\ref{eqn:iotacondition}) can be viewed as a notion of predictive stability in the sense of \citet{Wu:2005} on the tail conditional density and its derivatives. For many commonly used distributions in the modeling of tail dependence, the tail is usually heavier than the exponential and as a result the tail density function and its derivatives generally decay to zero faster than the tail probability. In such cases, condition (\ref{eqn:iotacondition}) typically does not require any additional assumption beyond what the tail adversarial stability requires in (A1). To illustrate this, we consider the moving-maximum process of \citet{Hall:Peng:Yao:2002} which was shown to be dense in the class of stationary processes whose finite-dimensional distributions are extreme-value of a given type; see also the discussion in \citet{Zhang:Zhang:Cui:2017}. In particular, for the quantile shifted moving-maximum process given in (\ref{eqn:einmm}), one can show that
\begin{equation*}
	\|\mathrm{E}\{\iota_{1,n}^{(q)}(\varepsilon; \mathcal{F}_{i}) \mid \mathcal{F}_{0}\} - \mathrm{E}\{\iota_{1,n}^{(q)}(\varepsilon; \mathcal{F}_{i}^\star) \mid \mathcal{F}_{0}^\star\}\| \leq {2^{1/2} a_{i+1}^{\kappa/2} \over a_0^q (Q_{e,1-\alpha_n}+\varepsilon)^{\kappa/2}} F_\epsilon^{(q)}\left({Q_{e,1-\alpha_n}+\varepsilon \over a_0}\right),
\end{equation*}
where $(Q_{e,1-\alpha_n}+\varepsilon)^{-\kappa} = O(\alpha_n)$ holds for any $|\varepsilon| \leq \varepsilon_0$ by the calculation in \citet{Hall:Peng:Yao:2002}. In addition, due to the algebraic decay of tail probabilities, it holds that $F_\epsilon^{(q)}\{a_0^{-1}(Q_{e,1-\alpha_n}+\varepsilon)\} = O[F_\epsilon^{(1)}\{a_0^{-1}(Q_{e,1-\alpha_n}+\varepsilon)\}]$ for $q = 1,\ldots,p$, where $F_\epsilon^{(1)}\{a_0^{-1}(Q_{e,1-\alpha_n}+\varepsilon)\} = O\{\alpha_n^{(\kappa+1)/\kappa}\}$ by the calculations in \citet{Hall:Peng:Yao:2002} and \citet{Zhang:2005}. As a result, if the coefficients satisfy $\sum_{l=0}^\infty a_l^{\kappa/2} < \infty$ as indicated by the tail adversarial stability in (A1), then
\begin{equation*}
	\sum_{i=0}^\infty \sup_{|\varepsilon| \leq \varepsilon_0}\|\mathrm{E}\{\iota_{1,n}^{(q)}(\varepsilon; \mathcal{F}_{i}) \mid \mathcal{F}_{0}\} - \mathrm{E}\{\iota_{1,n}^{(q)}(\varepsilon; \mathcal{F}_{i}^\star) \mid \mathcal{F}_{0}^\star\}\| = O\{\alpha_n^{1/2+(\kappa+1)/\kappa}\}
\end{equation*}
holds for $q = 1,\ldots,p$, and (\ref{eqn:iotacondition}) automatically follows as $\alpha_n^{1/2+(\kappa+1)/\kappa} = \alpha_n^{3/2+1/\kappa} = O(\alpha_n^{1/2})$. Let $\tilde \Sigma_n = n^{-1}\sum_{i=1}^n \gamma_{i,n}^{-1} \boldsymbol z_{i,n} \boldsymbol z_{i,n}^\top$, then Lemma \ref{lem:thetahatbound} provides the consistency of high quantile regression estimators when the auxiliary process is potentially heteroskedastic, which is needed as a preliminary for our main Bahadur representation result.

\begin{lemma}\label{lem:thetahatbound}
	Assume (A1)--(A4), $\alpha_n\to 0$ and $n\alpha_n\to \infty$. If
	\begin{equation*}
		\tau_n=(n\alpha_n)^{1/2}\frac{f_n(0)}{1-F_n(0)}\to \infty,
	\end{equation*}
	$\max_{1\leq i\leq n}|\boldsymbol z_{i,n}| = o(\tau_n)$,
	\begin{equation*}
		\max_{1\leq i\leq n}\sup_{|\boldsymbol \varphi|\leq c}\left|\frac{f_n(\gamma_{i,n}^{-1} \tau_n^{-1}\boldsymbol z_{i,n}^\top\boldsymbol \varphi)-f_n(0)}{f_n(0)}\right|\to 0
	\end{equation*}
	for any $c<\infty$, and eigenvalues of $\tilde \Sigma_n$ are bounded away from zero, then
	\begin{equation*}
		\hat{\boldsymbol \varphi}_{1-\alpha_n,n} - \boldsymbol \varphi_{1-\alpha_n,n} = O_p(\tau_n^{-1}).
	\end{equation*}
\end{lemma}


Using the preliminary results in Lemmas \ref{lem:Mn}--\ref{lem:thetahatbound}, we are now ready to establish the main Bahadur representation for high quantile regression estimators of tail dependent time series under possible heteroskedasticity in Theorem \ref{thm:Bahadurrepresentation}.

\begin{theorem}\label{thm:Bahadurrepresentation}
	Assume conditions of Lemmas \ref{lem:Mn}--\ref{lem:thetahatbound}. If $\sup_{u \in \mathbb R} f_{1,n}(u \mid \mathcal F_0) \leq c_{f,0}$ almost surely for some constant $c_{f,0} < \infty$ and there exists a constant $c_{f,1} < \infty$ such that
	\begin{equation}\label{eqn:fnprimecondition}
		\max_{1 \leq i \leq n} \sup_{|\boldsymbol \varphi| \leq c}\left|{f_n'(\gamma_{i,n}^{-1}\tau_n^{-1} \boldsymbol z_{i,n}^\top \boldsymbol \varphi) - f_n'(0) \over f_n'(0)}\right| \leq c_{f,1}
	\end{equation}
	holds for any $c < \infty$, then
	\begin{eqnarray}
		& & \tau_n \tilde{\Sigma}_n(\hat{\boldsymbol \varphi}_{1-\alpha_n,n}-\boldsymbol \varphi_{1-\alpha_n,n}) - (n\alpha_n)^{-1/2}\sum_{i=1}^n \psi_{1-\alpha_n}(\gamma_{i,n} e_{i,n}) \boldsymbol z_{i,n} \nonumber\\
		& = & O_p\left[{r_n^{1/2} \log n \over (n\alpha_n)^{1/4}} + {r_n \over \tau_n} \left\{1 + {f_n'(0) \over f_n(0)}\right\}\right].\label{eqn:Bahadurrepresentationorder}
	\end{eqnarray}
\end{theorem}

By Theorem \ref{thm:Bahadurrepresentation}, the nonlinear high quantile regression estimator $\hat{\boldsymbol \varphi}_{1-\alpha_n,n}$ can be well approximated by a suitably constructed linear form $(n\alpha_n)^{-1/2}\sum_{i=1}^n \psi_{1-\alpha_n}(\gamma_{i,n} e_{i,n}) \boldsymbol z_{i,n}$ with an explicit rate of approximation. Note that condition (\ref{eqn:fnprimecondition}) is similar in spirit to condition (9) used in Theorem 1 of \citet{Zhang:2021}, which essentially concerns the smoothness of the underlying tail distribution. In particular, since $\tau_n \to \infty$ as $n \to \infty$ and $\max_{1\leq i\leq n}|\boldsymbol z_{i,n}| = o(\tau_n)$, $\gamma_{i,n}^{-1}\tau_n^{-1} \boldsymbol z_{i,n}^\top \boldsymbol \varphi$ falls into a neighborhood of zero when $|\boldsymbol \varphi| \leq c$. As a result, $f_n'(\gamma_{i,n}^{-1}\tau_n^{-1} \boldsymbol z_{i,n}^\top \boldsymbol \varphi)$ is expected to be close to $f_n'(0)$ when the density function is continuously differentiable in the tail.

We shall here provide a discussion on the rate calculation in the high quantile Bahadur representation (\ref{eqn:Bahadurrepresentationorder}) and make a comparison with the fixed quantile setting. In particular, when the quantile level is fixed as in conventional quantile regression models, the quantile regression estimator and its approximated linear form in the traditional Bahadur representation are both scaled at the universal $n^{1/2}$ rate; see for example the results in \citet{Babu:Singh:1978}, \citet{Hesse:1990} and \citet{Wu:2005:Bahadur}. In comparison, the current growing quantile setting may alter the asymptotic behavior and makes the Bahadur representation follow a different rate. To be more specific, the constructed linear form in (\ref{eqn:Bahadurrepresentationorder}) now has a rate of $(n\alpha_n)^{1/2}$ while the high quantile regression estimator is scaled at a possibly different rate denoted by $\tau_n$ that can depend on an interplay between how extremal the quantile level is and the tail behavior of the underlying distribution. In addition, the explicit rate of approximation established in the high quantile Bahadur representation (\ref{eqn:Bahadurrepresentationorder}) involves two terms. The first term is at the rate of $(n\alpha_n)^{-1/4}$ up to a multiplicative logarithmic factor for applications \citep{Zhou:Shao:2013,Zhang:2021} with a bounded $r_n$. This term can be viewed as a high quantile analog of the optimal bound $n^{-1/4}$ in the conventional Bahadur representation when the quantile level is treated as fixed \citep{Bahadur:1966,Wu:2007:MEstimation}. The second term is at the rate of $\tau_n^{-1}$ under certain regularity conditions on the tail density and $r_n$. This term can be viewed as a higher-order term of the intrinsically nonlinear high quantile regression estimator. To be more specific, by Lemma \ref{lem:thetahatbound} the high quantile regression estimator satisfies $\hat{\boldsymbol \varphi}_{1-\alpha_n,n}-\boldsymbol \varphi_{1-\alpha_n,n} = O_p(\tau_n^{-1})$, and as a result if we multiply both sides of (\ref{eqn:Bahadurrepresentationorder}) by $\tau_n^{-1}$ then this second term becomes $\tau_n^{-2}$ which is a higher-order term of $\hat{\boldsymbol \varphi}_{1-\alpha_n,n}-\boldsymbol \varphi_{1-\alpha_n,n}$. We shall here use the moving-maximum process (\ref{eqn:einmm}) of \citet{Hall:Peng:Yao:2002} with different degrees of tail heaviness controlled by the shape parameter $\kappa$ to further illustrate the Bahadur representation result in Theorem \ref{thm:Bahadurrepresentation}. For simplicity, we consider the intercept model with $x_{i,n} \equiv 1$, then by construction $\Sigma_n = 1$ and $z_{i,n} = x_{i,n} = 1$. In this case, by properties of the moving-maximum process and the calculations in \citet{Hall:Peng:Yao:2002} and \citet{Zhang:2005}, one can show that
\begin{equation*}
	\tau_n = (n\alpha_n)^{1/2} \cdot {\kappa (1-\alpha_n) \over \alpha_n (\sum_{l=0}^\infty a_l^\kappa)^{1/\kappa}} \{-\log(1-\alpha_n)\}^{1+1/\kappa}
\end{equation*}
and
\begin{equation*}
	{f_n'(0) \over f_n(0)} = {\kappa \{-\log(1-\alpha_n)\}^{1+1/\kappa} - (\kappa+1)\{-\log(1-\alpha_n)\}^{1/\kappa} \over (\sum_{l=0}^\infty a_l^\kappa)^{1/\kappa}} = O(\alpha_n^{1/\kappa}).
\end{equation*}
Therefore, the high quantile Bahadur representation in (\ref{eqn:Bahadurrepresentationorder}) becomes
\begin{eqnarray*}
	& & \tau_n \tilde{\Sigma}_n(\hat{\boldsymbol \varphi}_{1-\alpha_n,n}-\boldsymbol \varphi_{1-\alpha_n,n}) - (n\alpha_n)^{-1/2}\sum_{i=1}^n \psi_{1-\alpha_n}(\gamma_{i,n} e_{i,n}) \boldsymbol z_{i,n} \\
	& = & O_p\left\{{\log n \over (n\alpha_n)^{1/4}} + {1 \over \tau_n}\right\} = O_p\left\{{\log n \over (n\alpha_n)^{1/4}} + {1 \over (n\alpha_n)^{1/2}\alpha_n^{1/\kappa}}\right\},
\end{eqnarray*}
where the first term $(n\alpha_n)^{-1/4}\log n$ dominates if $n^{1/4}\alpha_n^{1/4+1/\kappa}\log n \to \infty$ and otherwise the second term $(n\alpha_n)^{-1/2}\alpha_n^{-1/\kappa}$ dominates. In either case, the approximation bound is generally sharp enough for many statistical applications to reduce the original problem of studying the intrinsically nonlinear high quantile regression estimator to a more accessible problem about the linear form $(n\alpha_n)^{-1/2}\sum_{i=1}^n \psi_{1-\alpha_n}(\gamma_{i,n} e_{i,n}) \boldsymbol z_{i,n}$. Therefore, the high quantile Bahadur representation result developed in Theorem \ref{thm:Bahadurrepresentation} is expected to be transformative and useful in facilitating the study of many statistical problems regarding high quantile regression estimators. To illustrate, we in the following section consider an application to guide the development of a generative high quantile homogeneity test.

\section{A Generative High Quantile Homogeneity Test}\label{sec:application}
\subsection{Test Statistic and Asymptotic Theory}\label{subsec:applicationstatisticandtheory}
Suppose we are interested in comparing the regressor effect at two quantile levels $1-\alpha_{1,n}$ and $1-\alpha_{2,n}$, where both of them can grow to approach the boundary as $n \to \infty$. Let
\begin{equation*}
	\hat{\boldsymbol \beta}_{1-\alpha_{k,n},n} = \mathop{\mathrm{argmin}}_{\boldsymbol \beta \in \mathbb R^p} \sum_{i=1}^n \rho_{1-\alpha_{k,n}}(y_{i,n} - \boldsymbol x_{i,n}^\top \boldsymbol \beta),\quad k = 1,2,
\end{equation*}
be the associated high quantile regression estimators, and following the notation in (\ref{eqn:eintildedefinition}) we denote the associated auxiliary process by $\tilde e_{\alpha_{k,n},i,n} = \gamma_{\alpha_{k,n},i,n} e_{\alpha_{k,n},i,n}$, which can be potentially heteroskedastic with $\mathrm{pr}(e_{\alpha_{k,n},i,n} \leq 0) = 1-\alpha_{k,n}$. We shall here follow the same notation system and use $\boldsymbol \beta_{1-\alpha_{k,n},n}$, $\boldsymbol \varphi_{1-\alpha_{k,n},n}$, $\hat{\boldsymbol \varphi}_{1-\alpha_{k,n},n}$, $\tilde \Sigma_{\alpha_{k,n},n}$, $f_{\alpha_{k,n},n}(0)$, $f_{\alpha_{k,n},n}'(0)$, and $\tau_{\alpha_{k,n},n}$ to denote the corresponding quantities in Section \ref{sec:Bahadur} when $\alpha_n$ there is specifically taken as $\alpha_{k,n}$. Let
\begin{equation*}
\zeta_{i,n,k} = \{\alpha_{k,n}(1-\alpha_{k,n})\}^{-1/2}\psi_{1-\alpha_{k,n}}(\gamma_{\alpha_{k,n},i,n} e_{\alpha_{k,n},i,n})
\end{equation*}
and $\boldsymbol \zeta_{i,n} = (\zeta_{i,n,1},\zeta_{i,n,2})^\top$, then by the Bahadur representation in Theorem \ref{thm:Bahadurrepresentation} it suffices to focus on the linear form
\begin{equation*}
	\boldsymbol T_n = n^{-1/2} \sum_{i=1}^n \boldsymbol \zeta_{i,n} \otimes \boldsymbol z_{i,n},
\end{equation*}
where $\otimes$ denotes the Kronecker product. Write the pre-asymptotic covariance matrix $\Gamma_n = E(\boldsymbol T_n \boldsymbol T_n^\top)$, then Theorem \ref{thm:Tnclt} provides the multivariate central limit theorem for the tail linear quantity $\boldsymbol T_n$ with heteroskedastic auxiliary processes.

\begin{theorem}\label{thm:Tnclt}
	Assume that the conditions of Theorem \ref{lem:thetahatbound} hold for $\alpha_{k,n}$, $k = 1,2$. Then $\Gamma_n$ is positive semi-definite with bounded eigenvalues for all large n. If in addition the eigenvalues of $\Gamma_n$ are bounded away from zero for all large $n$, then
	\begin{equation*}
		\Gamma_n^{-1/2} \boldsymbol T_n \to_{d} N(0,\mathrm I_{2p \times 2p})
	\end{equation*}
	with $\mathrm I_{2p \times 2p}$ being the identity matrix in $\mathbb R^{2p \times 2p}$.
\end{theorem}

For any $l \in \mathbb Z$, let $l_+ = \max(l,0)$ and $l_- = \max(-l,0)$ be the positive and negative parts of $l$. Then, combining Theorem \ref{thm:Tnclt} with the Bahadur representation in Theorem \ref{thm:Bahadurrepresentation}, Corollary \ref{cor:varphinhatclt} provides the joint central limit theorem for high quantile regression estimators at different quantile levels.

\begin{corollary}\label{cor:varphinhatclt}
	Assume that the conditions of Theorem \ref{thm:Bahadurrepresentation} hold for $\alpha_{k,n}$ with
	\begin{equation*}
		{r_n^{1/2} \log n \over (n\alpha_{k,n})^{1/4}} + {r_n \over \tau_{\alpha_{k,n},n}} \left\{1 + {f_{\alpha_{k,n},n}'(0) \over f_{\alpha_{k,n},n}(0)}\right\} \to 0
	\end{equation*}
	for $k = 1,2$. If the limits
	\begin{equation*}
		\rho_{k_1,k_2}(l) = \lim_{n \to \infty} \mathrm{cor}(\mathbbm{1}_{\{e_{\alpha_{k_1,n},0,n} > 0\}}, \mathbbm{1}_{\{e_{\alpha_{k_2,n},l,n} > 0\}}),\quad \Upsilon(l) = \lim_{n \to \infty} n^{-1}\sum_{i=1+l_-}^{n-l_+} \boldsymbol z_{i,n} \boldsymbol z_{i+l,n}^\top
	\end{equation*}
	exist for each $l \in \mathbb Z$ and $k_1,k_2 \in \{1,2\}$, then the matrix
	\begin{equation*}
		\Gamma = \left\{\begin{array}{cc}
			\sum_{l \in \mathbb Z} \rho_{1,1}(l) \Upsilon(l) & \sum_{l \in \mathbb Z} \rho_{1,2}(l) \Upsilon(l) \\
			\sum_{l \in \mathbb Z} \rho_{2,1}(l) \Upsilon(l) & \sum_{l \in \mathbb Z} \rho_{2,2}(l) \Upsilon(l)
		\end{array}\right\}
	\end{equation*}
	is positive semi-definite with bounded eigenvalues for all large $n$. If in addition the eigenvalues of $\Gamma$ are bounded away from zero for all large $n$, then we have the central limit theorem
	\begin{equation*}
		\left\{\begin{array}{c}
			\tau_{\alpha_{1,n},n} \tilde \Sigma_{\alpha_{1,n},n}(\hat{\boldsymbol \varphi}_{1-\alpha_{1,n},n} - \boldsymbol \varphi_{1-\alpha_{1,n},n}) \\
			\tau_{\alpha_{2,n},n} \tilde \Sigma_{\alpha_{2,n},n}(\hat{\boldsymbol \varphi}_{1-\alpha_{2,n},n} - \boldsymbol \varphi_{1-\alpha_{2,n},n})
		\end{array}\right\} \to_d N(0,\Gamma).
	\end{equation*}
\end{corollary}

Note that $\Upsilon(-l) = \Upsilon(l)^\top$, $\rho_{k_1,k_2}(-l) = \rho_{k_1,k_2}(l)$ if $k_1 = k_2$ and $\rho_{k_1,k_2}(-l) = \rho_{k_2,k_1}(l)$ if $k_1 \neq k_2$, the matrix $\Gamma$ is automatically symmetric. Based on the joint central limit theorem provided in Corollary \ref{cor:varphinhatclt}, one can now construct a statistical test to determine the coefficient homogeneity at the two chosen high quantile levels $1-\alpha_{1,n}$ and $1-\alpha_{2,n}$ where both of them can grow to approach the upper boundary as $n \to \infty$. More generally, we consider testing the null hypothesis of
\begin{equation}\label{eqn:H0}
	H_0:\ A(\boldsymbol \beta_{1-\alpha_{1,n},n} - \boldsymbol \beta_{1-\alpha_{2,n},n}) = \boldsymbol c_A
\end{equation}
for some matrix $A \in \mathbb R^{p_0 \times p}$ and vector $\boldsymbol c_A \in \mathbb R^{p_0}$. For example, if we choose $A$ as the identity matrix and $\boldsymbol c_A$ as the zero vector, then (\ref{eqn:H0}) amounts to testing the coefficient homogeneity $\boldsymbol \beta_{1-\alpha_{1,n},n} = \boldsymbol \beta_{1-\alpha_{2,n},n}$ over the whole set of regressors. In practice, one may chose the matrix $A$ to focus on a certain set of regressors and change the vector $\boldsymbol c_A$ for different null considerations. Let
\begin{equation}\label{eqn:varthetanA}
	\boldsymbol \vartheta_n(A,\boldsymbol c_A) = A(\hat{\boldsymbol \beta}_{1-\alpha_{1,n},n} - \hat{\boldsymbol \beta}_{1-\alpha_{2,n},n}) - \boldsymbol c_A,
\end{equation}
Theorems \ref{cor:crossquantiletesting} and \ref{cor:crossquantiletestingdegenerate} provide the asymptotic distribution of the test statistic (\ref{eqn:varthetanA}) for situations when $\alpha_{1,n}$ and $\alpha_{2,n}$ are of the same rate and when they are of different rates, respectively.

\begin{theorem}\label{cor:crossquantiletesting}
	Assume conditions of Corollary \ref{cor:varphinhatclt}. If there exist $\tau_n \to \infty$, $0 < c_{\tau,1},c_{\tau,2} < \infty$ and matrices $\Sigma$, $\tilde \Sigma_1$ and $\tilde \Sigma_2$ with eigenvalues bounded away from zero and infinity such that
	\begin{equation*}
		\tau_{\alpha_{k,n},n}/\tau_n \to c_{\tau,k},\quad \Sigma_n \to \Sigma,\quad \tilde \Sigma_{\alpha_{k,n},n} \to \tilde \Sigma_k,
	\end{equation*}
	hold for $k = 1,2$, then under the null hypothesis (\ref{eqn:H0}),
	\begin{equation*}
		\tau_n^2 \boldsymbol \vartheta_n(A,\boldsymbol c_A)^\top (\tilde AB^{-1}\Gamma B^{-1}\tilde A^\top)^{-1} \boldsymbol \vartheta_n(A,\boldsymbol c_A) \to_d \chi^2_{p_0},
	\end{equation*}
	where
	\begin{equation*}
		\tilde A = (A,-A),\quad B = \left(\begin{array}{cc}
			c_{\tau,1} \tilde \Sigma_1 \Sigma^{1/2} & \\
			& c_{\tau,2} \tilde \Sigma_2 \Sigma^{1/2}
		\end{array}\right),
	\end{equation*}
	provided that the matrix $\tilde AB^{-1}\Gamma B^{-1}\tilde A^\top$ is invertible. In addition, the result continues to hold under the alternative if $\boldsymbol c_A$ is replaced by the actual $A(\boldsymbol \beta_{1-\alpha_{1,n},n} - \boldsymbol \beta_{1-\alpha_{2,n},n})$.
\end{theorem}

\begin{theorem}\label{cor:crossquantiletestingdegenerate}
	Assume conditions of Corollary \ref{cor:varphinhatclt}. If $\tau_{\alpha_{2,n},n}/\tau_{\alpha_{1,n},n} \to \infty$ and there exist matrices $\Sigma$, $\tilde \Sigma_1$ and $\tilde \Sigma_2$ with eigenvalues bounded away from zero and infinity such that
	\begin{equation*}
		\Sigma_n \to \Sigma,\quad \tilde \Sigma_{\alpha_{k,n},n} \to \tilde \Sigma_k,
	\end{equation*}
	hold for $k = 1,2$, then under the null hypothesis (\ref{eqn:H0}),
	\begin{equation*}
		\tau_{\alpha_{1,n},n}^2 \boldsymbol \vartheta_n(A,\boldsymbol c_A)^\top (A\Sigma^{-1/2} \tilde \Sigma_1^{-1} \tilde \Gamma \tilde \Sigma_1^{-1} \Sigma^{-1/2}A^\top)^{-1} \boldsymbol \vartheta_n(A,\boldsymbol c_A) \to_d \chi^2_{p_0},
	\end{equation*}
	where $\tilde \Gamma$ is the first $p \times p$ block of $\Gamma$ provided that the matrix $A\Sigma^{-1/2} \tilde \Sigma_1^{-1} \tilde \Gamma \tilde \Sigma_1^{-1} \Sigma^{-1/2}A^\top$ is invertible. In addition, the result continues to hold under the alternative if $\boldsymbol c_A$ is replaced by the actual $A(\boldsymbol \beta_{1-\alpha_{1,n},n} - \boldsymbol \beta_{1-\alpha_{2,n},n})$.
\end{theorem}

By Theorems \ref{cor:crossquantiletesting} and \ref{cor:crossquantiletestingdegenerate}, the asymptotic distribution of the test statistic (\ref{eqn:varthetanA}) can be different depending on whether $\alpha_{1,n}$ and $\alpha_{2,n}$ are of the same rate. Although this can be interesting from the theoretical perspective, it creates ambiguities for practitioners about which limiting distribution should be used in a given application. We shall in the following illustrate how the two cases can be unified together when testing coefficient homogeneity of regressors in the presence of an intercept for high quantile regression models.

\subsection{A Generative Approach for Unified Inference}\label{subsec:implementation}
Following the notation in Section \ref{subsec:applicationstatisticandtheory}, write $\boldsymbol x_{i,n} = (x_{i,n,1},\ldots,x_{i,n,p})^\top$, and we in this section consider testing the null model
\begin{equation}\label{eqn:applicationnullmodel}
	y_{i,n} = \beta_{1-\alpha_n,n,1} + \sum_{j=2}^p x_{i,n,j} \beta_{n,j} + e_{\alpha_n,i,n},\quad i = 1,\ldots,n,
\end{equation}
where $\beta_{1-\alpha_n,n,1}$ is the regression coefficient associated with the $(1-\alpha_n)$-th high quantile for the intercept term $x_{i,n,1} \equiv 1$ and $\beta_{n,2},\ldots,\beta_{n,p}$ are the corresponding regression coefficients for regressors $x_{i,n,2},\ldots,x_{i,n,p}$ which are assumed to be the same for different high quantiles. Therefore, under the null model (\ref{eqn:applicationnullmodel}), different high quantiles are associated with different intercepts but share the same coefficient effect from the regressors; see also the common-slope model in Section 2.1 of \citet{Wang:Li:He:2012}. We are interested in testing (\ref{eqn:applicationnullmodel}) against the alternative model
\begin{equation}\label{eqn:applicationalternativemodel}
	y_{i,n} = \beta_{1-\alpha_n,n,1} + \sum_{j=2}^p x_{i,n,j} \beta_{1-\alpha_n,n,j} + \tilde e_{\alpha_n,i,n},\quad i = 1,\ldots,n,
\end{equation}
where the regressor coefficients $\beta_{1-\alpha_n,n,2},\ldots,\beta_{1-\alpha_n,n,p}$ can also depend on $\alpha_n$ and thus be inhomogeneous for different high quantiles. For example, in the tropical cyclone data analysis of \citet{Elsner:Kossin:Jagger:2008}, the regressor is set as the time and it becomes desirable to test whether the increase of tropical cyclone intensities described by the slope is homogeneous at different quantile levels. Unlike the null model (\ref{eqn:applicationnullmodel}), the auxiliary process $(\tilde e_{\alpha_n,i,n})$ in the alternative model (\ref{eqn:applicationalternativemodel}) typically cannot be treated as a stationary array and it is more reasonable to consider the heteroskedastic case. To be more specific, by applying (\ref{eqn:applicationalternativemodel}) to quantile levels $1-\alpha_{1,n}$ and $1-\alpha_{2,n}$, the associated auxiliary processes then relate to each other by $\tilde e_{\alpha_{1,n},i,n} - \tilde e_{\alpha_{2,n},i,n} = \sum_{j=1}^p x_{i,n,j} (\beta_{1-\alpha_{2,n},n,j} - \beta_{1-\alpha_{1,n},n,j})$ which is not a constant and can be different for $i = 1,\ldots,n$. This can be accommodated by, for instance, setting $\tilde e_{\alpha_{k,n},i,n} = \gamma_{i,n} e_{\alpha_{k,n},i,n}$ for some $\gamma_{i,n}$ that is linear in $\boldsymbol x_{i,n}$; see also the location-scale shift quantile regression model in (2.14) of \citet{Wang:Li:He:2012}. The aforementioned article considered the problem of high quantile estimation separately under models (\ref{eqn:applicationnullmodel}) and (\ref{eqn:applicationalternativemodel}) for independent data with a polynomial-tail marginal distribution, and we shall here take advantage of the results developed in Sections \ref{sec:Bahadur} and \ref{subsec:applicationstatisticandtheory} to test the null model (\ref{eqn:applicationnullmodel}) against the alternative (\ref{eqn:applicationalternativemodel}) for a general class of tail dependent time series. The following algorithm provides a generative high quantile regression approach that unifies Theorems \ref{cor:crossquantiletesting} and \ref{cor:crossquantiletestingdegenerate} so that one does not need to distinguish whether $\alpha_{1,n}$ and $\alpha_{2,n}$ are of the same rate.
\begin{itemize}
	\item[(i)] Given high quantile levels $1-\alpha_{1,n}$ and $1-\alpha_{2,n}$, perform individual high quantile regression as in \citet{Zhang:2021} to obtain the coefficient estimates $\hat{\boldsymbol \beta}_{1-\alpha_{1,n},n}$ and $\hat{\boldsymbol \beta}_{1-\alpha_{2,n},n}$, along with their associated $\hat \tau_{\alpha_{1,n},n}$ and $\hat \tau_{\alpha_{2,n},n}$.
	\item[(ii)] Compute $\Upsilon_n(l) = n^{-1}\sum_{i=1+l_-}^{n-l_+} \boldsymbol z_{i,n} \boldsymbol z_{i+l,n}^\top$ and use the tapering estimate of \citet{Zhang:Xu:2024} at frequency zero to obtain estimates of $\Gamma_{k_1,k_2,n} = \sum_{l \in \mathbb Z} \rho_{k_1,k_2}(l) \Upsilon_n(l)$, denoted by $\hat \Gamma_{k_1,k_2,n}$, for $k_1,k_2 \in \{1,2\}$.
	\item[(iii)] Set $A = \{\boldsymbol 0_p,\mathrm I_{(p-1) \times (p-1)}\}$ as the block matrix to focus on the regression coefficients besides the intercept, and calculate
	\begin{equation*}
		\Phi = A\Sigma_n^{-1/2} \{\hat \tau_{\alpha_{1,n},n}^{-2}\hat \Gamma_{1,1,n} - \hat \tau_{\alpha_{1,n},n}^{-1} \hat \tau_{\alpha_{2,n},n}^{-1} (\hat \Gamma_{1,2,n} + \hat \Gamma_{2,1,n}) + \hat \tau_{\alpha_{2,n},n}^{-2}\hat \Gamma_{2,2,n}\} \Sigma_n^{-1/2}A^\top.
	\end{equation*}
	\item[(iv)] Generate independent samples $y_{i,n}^\circ$ from a given distribution and repeat (i)--(iii) on the simulation-assisted data $(\boldsymbol x_{i,n},y_{i,n}^\circ)$ to calculate the corresponding $\hat{\boldsymbol \beta}_{1-\alpha_{1,n},n}$, $\hat{\boldsymbol \beta}_{1-\alpha_{2,n},n}$, $\hat \tau_{\alpha_{1,n},n}$, $\hat \tau_{\alpha_{2,n},n}$, and $\Phi^\circ$.
	\item[(v)] Repeat (iv) and calculate the $(1-a)$-th quantile of
    \begin{equation*}
    2(\hat{\boldsymbol \beta}_{1-\alpha_{1,n},n}^\circ - \hat{\boldsymbol \beta}_{1-\alpha_{2,n},n}^\circ)^\top A^\top \{\Phi^\circ + (\Phi^\circ)^\top\}^{-1} A(\hat{\boldsymbol \beta}_{1-\alpha_{1,n},n}^\circ - \hat{\boldsymbol \beta}_{1-\alpha_{2,n},n}^\circ),
    \end{equation*}
    denoted by $Q_{1-a}^\circ$.
	\item[(vi)] Reject the null model (\ref{eqn:applicationnullmodel}) at level $1-a$ if
	\begin{equation*}
		2(\hat{\boldsymbol \beta}_{1-\alpha_{1,n},n} - \hat{\boldsymbol \beta}_{1-\alpha_{2,n},n})^\top A^\top (\Phi + \Phi^\top)^{-1} A(\hat{\boldsymbol \beta}_{1-\alpha_{1,n},n} - \hat{\boldsymbol \beta}_{1-\alpha_{2,n},n}) > Q_{1-a}^\circ.
	\end{equation*}
\end{itemize}

In addition to unifying Theorems \ref{cor:crossquantiletesting} and \ref{cor:crossquantiletestingdegenerate}, the above generative algorithm can also lead to improved finite-sample performance as it mimics how the uncertainty is accumulated from different sources in the target statistic. Since the effect of dependence is pivotalized in the asymptotic distribution provided in Theorems \ref{cor:crossquantiletesting} and \ref{cor:crossquantiletestingdegenerate}, the mechanism here is in connection with the finite-sample approximation of self-normalized distributions \citep{Shao:2010,Zhang:Lavitas:2018}. In practice, the distribution from which the simulation-assisted samples are generated can be general as the resulting asymptotic distribution will remain the same by Theorems \ref{cor:crossquantiletesting} and \ref{cor:crossquantiletestingdegenerate}. As a rule of thumb, we recommend generating the simulation-assisted samples from the Fr\'{e}chet family or the generalized Pareto distribution with an estimated parameter \citep{deZeaBermudez:Kotz:2010:PartI,deZeaBermudez:Kotz:2010:PartII,Cao:Gao:Shao:Sriram:Wang:Wen:Zhang:2025}. It can be seen from our numerical experiments in Sections \ref{subsec:simulation} and \ref{subsec:dataanalysisNASDAQ} that the proposed generative high quantile homogeneity test seems to perform reasonably well.

\subsection{A Simulation Study}\label{subsec:simulation}
We in this section provide a Monte Carlo simulation study to examine the finite-sample performance of the generative high quantile regression approach proposed in Section \ref{subsec:implementation} for testing coefficient homogeneity of regressors at high quantile levels $1-\alpha_{1,n}$ and $1-\alpha_{2,n}$. For this, we consider the trend analysis setting of \citet{Elsner:Kossin:Jagger:2008} and \citet{Rhines:McKinnon:Tingley:Huybers:2017} in which $x_{i,n,j} = t_{i,n}^{j-1}$, $j = 1,\ldots,p$, are taken as polynomials of the time $t_{i,n} = i/n$, $i = 1,\ldots,n$, and model (\ref{eqn:applicationnullmodel}) becomes
\begin{equation}\label{eqn:modelsimulation}
	y_{i,n} = \beta_{1-\alpha_{k,n},n,1} + \sum_{j=2}^p \beta_{n,j} t_{i,n}^{j-1} + e_{\alpha_{k,n},i,n},\quad i = 1,\ldots,n.
\end{equation}
To examine the performance under nonnegligible tail dependence, we follow \citet{Zhang:2021} and consider the moving-maximum \citep{Hall:Peng:Yao:2002} auxiliary process
\begin{equation*}
	e_{\alpha_{k,n},i,n} = \max(\epsilon_i,a\epsilon_{i-1}) - b_{n,k},
\end{equation*}
where $(\epsilon_i)$ is a sequence of independent exponential random variables and $b_{n,k}$ is the shift to properly center the auxiliary process so that $\mathrm{pr}(e_{\alpha_{k,n},i,n} > 0) = \alpha_{k,n}$ for $k = 1,2$; see also \citet{Koenker:Bassett:1978}. The case when $a = 0$ relates to independent data, and a nonzero choice will lead to tail dependent data. Let $n \in \{1000,2000,5000\}$, $a \in \{0,1\}$ and $\beta_{n,j} = j$, for each realization, we apply the generative algorithm proposed in Section \ref{subsec:implementation} to test if the high quantile trends at $1-\alpha_{1,n} = 0.95$ and $1-\alpha_{2,n} = 0.99$ follow the same polynomial pattern except for an intercept shift, namely if the change in high quantiles with respect to time is homogeneous at high quantile levels 95\% and 99\%. We also report results based on the theoretical $\chi^2$-distributions established in Theorems \ref{cor:crossquantiletesting} and \ref{cor:crossquantiletestingdegenerate} for a comparison. In practice, it can be ambiguous to determine which $\chi^2$-distribution from Theorems \ref{cor:crossquantiletesting} and \ref{cor:crossquantiletestingdegenerate} should be used to better represent the underlying asymptotic setting, while the generative algorithm described in Section \ref{subsec:implementation} provides a unified inference procedure with potential finite-sample improvement. It can be seen from Table \ref{tab:simulation9599} that noticeable size distortions can be observed for the theoretical $\chi^2$-approximations, especially when the sample size is small and when multiple components are being tested. For example, when $n = 1000$ and $p = 4$, the empirical coverage probabilities for model (\ref{eqn:modelsimulation}) with $a = 1$ are $(0.350,0.434,0.581)$ and $(0.363,0.443,0.595)$ at the $(90\%,95\%,99\%)$ nominal levels for the $\chi^2$-approximations using Theorems \ref{cor:crossquantiletesting} and \ref{cor:crossquantiletestingdegenerate} respectively. Note that when $n = 1000$ there are only about 10 data points in the tail of the 99\% high quantile, which can be considered as a very difficult scenario. If we increase the sample size, however, the amount of size distortions generally becomes smaller. In particular, the empirical coverage probabilities of the aforementioned difficult scenario can improve from $(0.350,0.434,0.581)$ and $(0.363,0.443,0.595)$ to $(0.676,0.760,0.878)$ and $(0.743,0.817,0.914)$ when $n$ increases from $1000$ to $5000$. On the other hand, the generative procedure in Section \ref{subsec:implementation} seems to be able to deliver an improved finite-sample performance especially when the sample size is relatively small, and the corresponding empirical coverage probabilities are generally close to their nominal levels. We in Table \ref{tab:simulation9099} also consider the situation when $1-\alpha_{1,n} = 0.9$ and $1-\alpha_{2,n} = 0.99$, and the aforementioned simulation findings seem to continue to hold.

{\renewcommand{\arraystretch}{0.7}
\begin{table}[!t]\centering
	\caption{\small Empirical coverage probabilities of the generative procedure described in Section \ref{subsec:implementation} for testing coefficient homogeneity at high quantile levels 95\% and 99\%. Results based on the theoretical $\chi^2$-approximations of Theorems \ref{cor:crossquantiletesting} and \ref{cor:crossquantiletestingdegenerate} are also reported and denoted by theoretical $\chi^2$ (I) and theoretical $\chi^2$ (II) respectively.}\label{tab:simulation9599}
	\begin{tabular}{llcccccccccccc}
		\hline
		& & & \multicolumn{3}{c}{generative algorithm} & & \multicolumn{3}{c}{theoretical $\chi^2$ (I)} & & \multicolumn{3}{c}{theoretical $\chi^2$ (II)} \\ \cline{4-6} \cline{8-10} \cline{12-14}
		$n$ & $p$ & & 90\% & 95\% & 99\% & & 90\% & 95\% & 99\% & & 90\% & 95\% & 99\% \\
		\hline
		\vspace{-1em}\\
		& & & \multicolumn{11}{c}{$a = 1$} \\
		1000 & 2 & & 0.868 & 0.926 & 0.982 & & 0.872 & 0.904 & 0.949 & & 0.881 & 0.914 & 0.954 \\
		& 3 & & 0.848 & 0.913 & 0.978 & & 0.627 & 0.698 & 0.810 & & 0.645 & 0.718 & 0.821 \\
		& 4 & & 0.817 & 0.893 & 0.972 & & 0.350 & 0.434 & 0.581 & & 0.363 & 0.443 & 0.595 \\
		2000 & 2 & & 0.880 & 0.936 & 0.985 & & 0.921 & 0.947 & 0.976 & & 0.934 & 0.956 & 0.981 \\
		& 3 & & 0.868 & 0.928 & 0.982 & & 0.745 & 0.809 & 0.894 & & 0.780 & 0.837 & 0.910 \\
		& 4 & & 0.853 & 0.920 & 0.984 & & 0.503 & 0.595 & 0.742 & & 0.548 & 0.639 & 0.775 \\
		5000 & 2 & & 0.897 & 0.940 & 0.986 & & 0.951 & 0.970 & 0.989 & & 0.965 & 0.980 & 0.993 \\
		& 3 & & 0.881 & 0.938 & 0.985 & & 0.842 & 0.894 & 0.954 & & 0.881 & 0.923 & 0.970 \\
		& 4 & & 0.874 & 0.934 & 0.983 & & 0.676 & 0.760 & 0.878 & & 0.743 & 0.817 & 0.914 \\
		\vspace{-1em}\\
		& & & \multicolumn{11}{c}{$a = 0$} \\
		1000 & 2 & & 0.925 & 0.967 & 0.996 & & 0.924 & 0.950 & 0.978 & & 0.937 & 0.959 & 0.982 \\
		& 3 & & 0.939 & 0.971 & 0.997 & & 0.764 & 0.830 & 0.912 & & 0.794 & 0.852 & 0.922 \\
		& 4 & & 0.948 & 0.978 & 0.998 & & 0.527 & 0.617 & 0.764 & & 0.566 & 0.656 & 0.793 \\
		2000 & 2 & & 0.920 & 0.962 & 0.993 & & 0.947 & 0.966 & 0.988 & & 0.959 & 0.976 & 0.991 \\
		& 3 & & 0.938 & 0.972 & 0.996 & & 0.832 & 0.886 & 0.947 & & 0.866 & 0.910 & 0.958 \\
		& 4 & & 0.944 & 0.979 & 0.998 & & 0.640 & 0.725 & 0.849 & & 0.696 & 0.775 & 0.884 \\
		5000 & 2 & & 0.921 & 0.961 & 0.994 & & 0.971 & 0.983 & 0.995 & & 0.982 & 0.990 & 0.998 \\
		& 3 & & 0.932 & 0.970 & 0.994 & & 0.894 & 0.936 & 0.976 & & 0.929 & 0.958 & 0.985 \\
		& 4 & & 0.930 & 0.967 & 0.994 & & 0.755 & 0.832 & 0.924 & & 0.820 & 0.881 & 0.953 \\
		\vspace{-1em}\\
		\hline
	\end{tabular}
\end{table}}

{\renewcommand{\arraystretch}{0.7}
\begin{table}[!t]\centering
	\caption{\small Empirical coverage probabilities of the generative procedure described in Section \ref{subsec:implementation} for testing coefficient homogeneity at high quantile levels 90\% and 99\%. Results based on the theoretical $\chi^2$-approximations of Theorems \ref{cor:crossquantiletesting} and \ref{cor:crossquantiletestingdegenerate} are also reported and denoted by theoretical $\chi^2$ (I) and theoretical $\chi^2$ (II) respectively.}\label{tab:simulation9099}
	\begin{tabular}{llcccccccccccc}
		\hline
		& & & \multicolumn{3}{c}{generative algorithm} & & \multicolumn{3}{c}{theoretical $\chi^2$ (I)} & & \multicolumn{3}{c}{theoretical $\chi^2$ (II)} \\ \cline{4-6} \cline{8-10} \cline{12-14}
		$n$ & $p$ & & 90\% & 95\% & 99\% & & 90\% & 95\% & 99\% & & 90\% & 95\% & 99\% \\
		\hline
		\vspace{-1em}\\
		& & & \multicolumn{11}{c}{$a = 1$} \\
		1000 & 2 & & 0.874 & 0.932 & 0.983 & & 0.869 & 0.903 & 0.949 & & 0.868 & 0.901 & 0.945 \\
		& 3 & & 0.862 & 0.926 & 0.982 & & 0.617 & 0.694 & 0.806 & & 0.606 & 0.683 & 0.794 \\
		& 4 & & 0.845 & 0.920 & 0.983 & & 0.329 & 0.412 & 0.567 & & 0.314 & 0.393 & 0.540 \\
		2000 & 2 & & 0.883 & 0.937 & 0.984 & & 0.918 & 0.945 & 0.974 & & 0.923 & 0.948 & 0.975 \\
		& 3 & & 0.874 & 0.930 & 0.985 & & 0.738 & 0.804 & 0.889 & & 0.748 & 0.811 & 0.894 \\
		& 4 & & 0.852 & 0.916 & 0.985 & & 0.479 & 0.571 & 0.722 & & 0.489 & 0.582 & 0.731 \\
		5000 & 2 & & 0.892 & 0.947 & 0.987 & & 0.952 & 0.970 & 0.990 & & 0.959 & 0.974 & 0.991 \\
		& 3 & & 0.881 & 0.938 & 0.985 & & 0.841 & 0.891 & 0.952 & & 0.858 & 0.903 & 0.960 \\
		& 4 & & 0.873 & 0.930 & 0.985 & & 0.657 & 0.747 & 0.873 & & 0.687 & 0.774 & 0.887 \\
		\vspace{-1em}\\
		& & & \multicolumn{11}{c}{$a = 0$} \\
		1000 & 2 & & 0.926 & 0.964 & 0.996 & & 0.921 & 0.946 & 0.975 & & 0.926 & 0.950 & 0.977 \\
		& 3 & & 0.946 & 0.977 & 0.998 & & 0.749 & 0.813 & 0.900 & & 0.756 & 0.821 & 0.904 \\
		& 4 & & 0.949 & 0.980 & 0.998 & & 0.504 & 0.598 & 0.743 & & 0.514 & 0.608 & 0.746 \\
		2000 & 2 & & 0.920 & 0.964 & 0.994 & & 0.947 & 0.965 & 0.987 & & 0.953 & 0.970 & 0.989 \\
		& 3 & & 0.937 & 0.974 & 0.996 & & 0.824 & 0.879 & 0.944 & & 0.839 & 0.890 & 0.950 \\
		& 4 & & 0.946 & 0.978 & 0.998 & & 0.630 & 0.720 & 0.846 & & 0.654 & 0.740 & 0.858 \\
		5000 & 2 & & 0.910 & 0.960 & 0.994 & & 0.967 & 0.981 & 0.994 & & 0.974 & 0.984 & 0.995 \\
		& 3 & & 0.928 & 0.967 & 0.994 & & 0.889 & 0.930 & 0.974 & & 0.906 & 0.941 & 0.980 \\
		& 4 & & 0.937 & 0.971 & 0.997 & & 0.759 & 0.837 & 0.928 & & 0.789 & 0.861 & 0.939 \\
		\vspace{-1em}\\
		\hline
	\end{tabular}
\end{table}}

\subsection{A Financial Data Application}\label{subsec:dataanalysisNASDAQ}
\begin{figure}[!t]\centering
	\begin{center}
		\includegraphics[width=\textwidth]{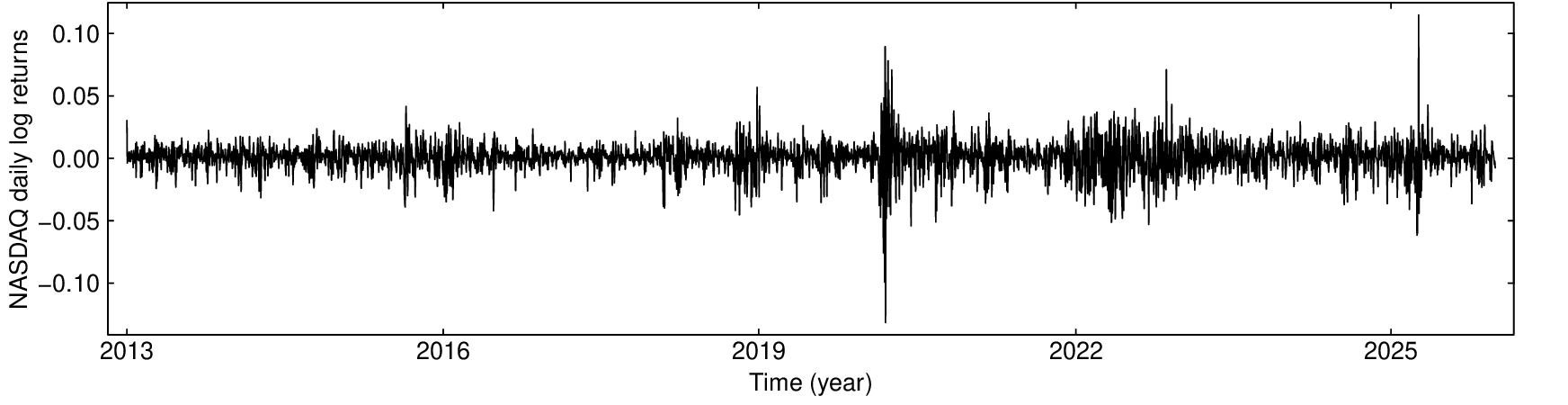}
		\caption{\small Daily log returns calculated from the NASDAQ composite index value at market close for the period of 2013--2025.}\label{fig:NASDAQ}
	\end{center}
\end{figure}

We in this section further illustrate the developed results by a financial data that contains the daily NASDAQ composite index value at market close for the period of 2013--2025. The data is available from the Federal Reserve Bank of St. Louis at \url{https://fred.stlouisfed.org/series/NASDAQCOM}, and a time series plot of the log return series is provided in Figure \ref{fig:NASDAQ}. For daily stock returns, central tendency measures such as the mean and median are often close to zero, making simple trend specifications such as a constant or linear form particularly appealing. As suggested in the \citet{IMF:2025} annual report, financial risks have increased due to tightening conditions, trade uncertainty, and vulnerabilities in capital markets, among other factors. We shall here consider the problem of testing if the increase captured by the slope of a linear trend is homogeneous in different quantiles, especially among lower tail quantiles that relate to extremal loss risks. Note that the lower tail quantile relates to the Value at Risk (VaR), an important risk metric, and can be handled by focusing on the upper tail quantile of the negative transformed time series. We model the data as a general tail dependent process, and the generative testing procedure in Section \ref{subsec:implementation} yields a $p$-value of 0.135 for a comparison between the 1\% and 5\% tail quantiles. As a result, the largest 1\% and 5\% losses seem to exhibit similar increasing patterns over time, differing primarily by a vertical shift that accommodate their different baselines. If we apply the test to compare the 5\% tail quantile with the median, however, the $p$-value then becomes 0.001 with the tail quantile exhibiting a sharper increase over time than the median. Our analysis reveals that the increase in financial risk reflected by stock returns is not uniform across the distribution. Lower tail quantiles corresponding to worst-case losses seem to display a significantly steeper upward trend than the median, implying that extreme downside risk is worsening at a faster pace than risks reflected by central tendency measures. In another words, the worst case is getting worse faster. This has important implications for financial institutions and the insurance industry to proactively adapt their risk management frameworks to better mitigate escalating tail risks.

\section{Discussion}\label{sec:discussion}
In this article, we establish the Bahadur representation of high quantile regression estimators for a general class of tail dependent time series. The result in Section \ref{sec:Bahadur} states that, after proper centering and scaling, the intrinsically nonlinear high quantile regression estimators can be well approximated by suitably constructed linear forms of tail dependent processes with an explicit rate of approximation. Compared with the result in \citet{Zhang:2021}, a much more sophisticated theoretical treatment involving the control of tail probabilities of martingales and a deep probabilistic chaining argument is needed to obtain the uniform negligibility with an explicit bound over an uncountable set that was otherwise not needed in \citet{Zhang:2021}. Since linear forms are generally much easier to deal with, the Bahadur representation is expected to be useful in facilitating the study of many statistical problems regarding high quantile regression estimators. As an illustration, we in Section \ref{sec:application} consider an application to guide the development of a generative high quantile homogeneity test to determine the statistical significance between high quantile estimators obtained from the same data but at different high quantile levels. Taking advantage of the newly developed Bahadur representation, the problem is reduced to the joint asymptotic analysis of linear forms, and we establish a dichotomous asymptotic distribution of the difference statistic depending on the relative extremeness of the two given high quantile levels. A testing procedure using generative samples is then proposed to unify the two cases, and numerical experiments including Monte Carlo simulations and a real data example are further provided to illustrate the developed results. It is expected that the Bahadur representation can lead to meaningful developments for other problems involving high quantile regression estimators, and we shall leave such research topics to future explorers.


\section*{Supplementary Material}
Supplementary material contains technical proofs of results in Sections \ref{sec:Bahadur} and \ref{sec:application}.

\section*{Disclosure statement}
The authors report there are no competing interests to declare.

\setlength{\bibsep}{0mm}
\bibliography{BibliographyTing,BibliographyTingBRHQ}

@ARTICLE{Bahadur:1966,
    author    = {R.~R. Bahadur},
    title     = {A note on quantiles in large samples},
    journal   = {The Annals of Mathematical Statistics},
    volume    = {37},
    number    = {3},
    pages     = {577--580},
    year      = {1966},
    month     = {June},
}

@ARTICLE{Bai:Zhang:2024,
    author    = {Shuyang Bai and Ting Zhang},
    title     = {Tail adversarial stability for regularly varying linear processes and their extensions},
    journal   = {Extremes},
    volume    = {27},
    number    = {1},
    pages     = {33--65},
    year      = {2024},
    month     = {March},
}

@ARTICLE{Bai:Rao:Wu:1992,
    author    = {Z.~D. Bai and C.~Radhakrishna Rao and Y. Wu},
    title     = {M-estimation of multivariate linear regression parameters under a convex discrepancy function},
    journal   = {Statistica Sinica},
    volume    = {2},
    number    = {1},
    pages     = {237--254},
    year      = {1992},
    month     = {January},
}

@ARTICLE{deZeaBermudez:Kotz:2010:PartI,
    author    = {P. {de Zea Bermudez} and Samuel Kotz},
    title     = {Parameter estimation of the generalized Pareto distribution-Part {I}},
    journal   = {Journal of Statistical Planning and Inference},
    volume    = {140},
    issue     = {6},
    pages     = {1353--1373},
    year      = {2010},
    month     = {June},
}

@ARTICLE{deZeaBermudez:Kotz:2010:PartII,
    author    = {P. {de Zea Bermudez} and Samuel Kotz},
    title     = {Parameter estimation of the generalized Pareto distribution-Part {II}},
    journal   = {Journal of Statistical Planning and Inference},
    volume    = {140},
    issue     = {6},
    pages     = {1374--1388},
    year      = {2010},
    month     = {June},
}

@ARTICLE{Cao:Gao:Shao:Sriram:Wang:Wen:Zhang:2025,
    author    = {Hanyue Cao and Jingying Gao and Yu Shao and T.~N. Sriram and Weiliang Wang and Fei Wen and Ting Zhang},
    title     = {Tail index estimation for tail adversarial stable time series with an application to high-dimensional tail clustering},
    journal   = {Journal of Time Series Analysis},
    volume    = {},
    number    = {},
    pages     = {forthcoming},
    year      = {2025},
}

@ARTICLE{Chernozhukov:2005,
    author    = {Victor Chernozhukov},
    title     = {Extremal quantile regression},
    journal   = {The Annals of Statistics},
    volume    = {33},
    number    = {2},
    pages     = {806--839},
    year      = {2005},
}

@ARTICLE{Chernozhukov:FernandezVal:2011,
    author    = {Victor Chernozhukov and Iv\'{a}n Fern\'{a}ndez-Val},
    title     = {Inference for extremal conditional quantile models, with an application to market and birthweight risks},
    journal   = {The Review of Economic Studies},
    volume    = {78},
    issue     = {2},
    pages     = {559--589},
    year      = {2011},
    month     = {April},
}

@ARTICLE{Daouia:Stupfler:Usseglio-Carleve:2023,
    author    = {Abdelaati Daouia and Gilles Claude Stupfler and Antoine Usseglio-Carleve},
    title     = {Inference for extremal regression with dependent heavy-tailed data},
    journal   = {The Annals of Statistics},
    volume    = {51},
    number    = {5},
    pages     = {2040--2066},
    year      = {2023},
    month     = {October},
}

@ARTICLE{Elsner:Kossin:Jagger:2008,
    author    = {James B. Elsner and James P. Kossin and Thomas H. Jagger},
    title     = {The increasing intensity of the strongest tropical cyclones},
    journal   = {Nature},
    volume    = {455},
    pages     = {92--95},
    year      = {2008},
    month     = {September},
}

@ARTICLE{Freedman:1975,
    author    = {David A. Freedman},
    title     = {On tail probabilities for martingales},
    journal   = {The Annals of Probability},
    volume    = {3},
    issue     = {1},
    pages     = {100--118},
    year      = {1975},
    month     = {February},
}

@ARTICLE{Hall:Peng:Yao:2002,
    author    = {Peter Hall and Liang Peng and Qiwei Yao},
    title     = {Moving-maximum models for extrema of time series},
    journal   = {Journal of Statistical Planning and Inference},
    volume    = {103},
    issue     = {1--2},
    pages     = {51--63},
    year      = {2002},
    month     = {April},
}

@ARTICLE{Heffernan:Tawn:Zhang:2007,
    author    = {Janet E. Heffernan and Jonathan A. Tawn and Zhengjun Zhang},
    title     = {Asymptotically (in)dependent multivariate maxima of moving maxima processes},
    journal   = {Extremes},
    volume    = {10},
    issue     = {2},
    pages     = {57--82},
    year      = {2007},
    month     = {June},
}

@ARTICLE{Hesse:1990,
    author    = {C.~H. Hesse},
    title     = {A Bahadur-type representation for empirical quantiles of a large class of stationary, possibly infinite-variance, linear processes},
    journal   = {The Annals of Statistics},
    volume    = {18},
    number    = {3},
    pages     = {1188--1202},
    year      = {1990},
    month     = {September},
}

@ARTICLE{Ho:Hsing:1996,
    author    = {Hwai-Chung Ho and Tailen Hsing},
    title     = {On the asymptotic expansion of the empirical process of long-memory moving averages},
    journal   = {The Annals of Statistics},
    volume    = {24},
    number    = {3},
    pages     = {992--1024},
    year      = {1996},
    month     = {June},
}

@BOOK{Koenker:2005,
    author    = {Roger Koenker},
    title     = {Quantile Regression},
    publisher = {Cambridge University Press},
    year      = {2005},
    address   = {Cambridge},
}

@ARTICLE{Koenker:Bassett:1978,
    author    = {Roger Koenker and Gilbert Jr. Bassett},
    title     = {Regression quantiles},
    journal   = {Econometrica},
    volume    = {46},
    number    = {1},
    pages     = {33--50},
    year      = {1978},
    month     = {January}
}

@ARTICLE{Koenker:Bassett:1982,
    author    = {Roger Koenker and Gilbert Jr. Bassett},
    title     = {Robust tests for heteroscedasticity based on regression quantiles},
    journal   = {Econometrica},
    volume    = {50},
    number    = {1},
    pages     = {43--61},
    year      = {1982},
    month     = {January}
}

@ARTICLE{Rhines:McKinnon:Tingley:Huybers:2017,
    author    = {Andrew Rhines and Karen A. McKinnon and Martin P. Tingley and Peter Huybers},
    title     = {Seasonally resolved distributional trends of North American temperatures show contraction of winter variability},
    journal   = {Journal of Climate},
    volume    = {30},
    number    = {3},
    pages     = {1139--1157},
    year      = {2017},
    month     = {February},
}

@ARTICLE{Rosenblatt:1956,
    author    = {Murray Rosenblatt},
    title     = {A central limit theorem and a strong mixing condition},
    journal   = {Proceedings of the National Academy of Sciences of the United States of America},
    volume    = {42},
    number    = {1},
    pages     = {43--47},
    year      = {1956},
    month     = {January},
}

@ARTICLE{Shao:2010,
    author    = {Xiaofeng Shao},
    title     = {A self-normalized approach to confidence interval construction in time series},
    journal   = {Journal of the Royal Statistical Society: Series B (Statistical Methodology)},
    volume    = {72},
    issue     = {3},
    pages     = {343--366},
    year      = {2010},
}

@BOOK{Tong:1990,
    author    = {Howell Tong},
    title     = {Non-Linear Time Series: A Dynamical System Approach (Oxford Statistical Science Series, 6)},
    publisher = {Oxford University Press},
    year      = {1990},
    address   = {U.K.},
}

@ARTICLE{Wang:Li:He:2012,
    author    = {Huixia Judy Wang and Deyuan Li and Xuming He},
    title     = {Estimation of high conditional quantiles for heavy-tailed distributions},
    journal   = {Journal of the American Statistical Association},
    volume    = {107},
    number    = {500},
    pages     = {1453--1464},
    year      = {2012},
    month     = {December},
}

@BOOK{Wiener:1958,
    author    = {Norbert Wiener},
    title     = {Nonlinear Problems In Random Theory (Technology Press Research Monographs)},
    publisher = {The MIT Press},
    year      = {1958},
    address   = {Cambridge},
}

@ARTICLE{Wu:2005,
    author    = {Wei Biao Wu},
    title     = {Nonlinear system theory: another look at dependence},
    journal   = {Proceedings of the National Academy of Sciences of the United States of America},
    volume    = {102},
    number    = {40},
    pages     = {14150--14154},
    year      = {2005},
    month     = {October},
}

@ARTICLE{Wu:2005:Bahadur,
    author    = {Wei Biao Wu},
    title     = {On the {B}ahadur representation of sample quantiles for dependent sequences},
    journal   = {The Annals of Statistics},
    volume    = {33},
    number    = {4},
    pages     = {1934--1963},
    year      = {2005},
    month     = {August},
}

@ARTICLE{Wu:2007:MEstimation,
    author    = {Wei Biao Wu},
    title     = {M-estimation of linear models with dependent errors},
    journal   = {The Annals of Statistics},
    volume    = {35},
    number    = {2},
    pages     = {495--521},
    year      = {2007},
}

@ARTICLE{Wu:2011,
    author    = {Wei Biao Wu},
    title     = {Asymptotic theory for stationary processes},
    journal   = {Statistics and Its Interface},
    volume    = {4},
    number    = {2},
    pages     = {207--226},
    year      = {2011},
}

@ARTICLE{Wu:Zhou:2017,
    author    = {Weichi Wu and Zhou Zhou},
    title     = {Nonparametric inference for time-varying coefficient quantile regression},
    journal   = {Journal of Business \& Economic Statistics},
    volume    = {35},
    issue     = {1},
    pages     = {98--109},
    year      = {2017},
}

@ARTICLE{Zhang:2021,
    author    = {Ting Zhang},
    title     = {High-quantile regression for tail-dependent time series},
    journal   = {Biometrika},
    volume    = {108},
    number    = {1},
    pages     = {113--126},
    year      = {2021},
}

@ARTICLE{Zhang:2022,
    author    = {Ting Zhang},
    title     = {Asymptotics of sample tail autocorrelations for tail dependent time series: phase transition and visualization},
    journal   = {Biometrika},
    volume    = {109},
    number    = {2},
    pages     = {521--534},
    year      = {2022},
}

@ARTICLE{Zhang:Lavitas:2018,
    author    = {Ting Zhang and Liliya Lavitas},
    title     = {Unsupervised self-normalized change-point testing for time series},
    journal   = {Journal of the American Statistical Association},
    volume    = {113},
    issue     = {522},
    pages     = {637--648},
    year      = {2018},
}

@ARTICLE{Zhang:Xu:2024,
    author    = {Ting Zhang and Beibei Xu},
    title     = {Tail spectral density estimation and its uncertainty quantification: another look at tail dependent time series analysis},
    journal   = {Journal of the American Statistical Association},
    volume    = {119},
    issue     = {546},
    pages     = {1424--1433},
    year      = {2024},
}

@ARTICLE{Zhang:2005,
    author    = {Zhengjun Zhang},
    title     = {A new class of tail-dependent time series models and its applications in financial time series},
    journal   = {Advances in Econometrics},
    volume    = {20},
    pages     = {323--358},
    year      = {2005},
}

@ARTICLE{Zhang:Smith:2004,
    author    = {Zhengjun Zhang and Richard L. Smith},
    title     = {The behavior of multivariate maxima of moving maxima processes},
    journal   = {Journal of Applied Probability},
    volume    = {41},
    issue     = {4},
    pages     = {1113--1123},
    year      = {2004},
    month     = {December},
}

@ARTICLE{Zhang:Zhang:Cui:2017,
    author    = {Zhengjun Zhang and Chunming Zhang and Qiurong Cui},
    title     = {Random threshold driven tail dependence measures with application to precipitation data analysis},
    journal   = {Statistica Sinica},
    volume    = {27},
    issue     = {2},
    pages     = {685--709},
    year      = {2017},
    month     = {April},
}

@ARTICLE{Zhou:Shao:2013,
    author    = {Zhou Zhou and Xiaofeng Shao},
    title     = {Inference for linear models with dependent errors},
    journal   = {Journal of the Royal Statistical Society: Series B (Statistical Methodology)},
    volume    = {75},
    issue     = {2},
    pages     = {323--343},
    year      = {2013},
}

@ARTICLE{Babu:Singh:1978,
    author    = {Gutti Jogesh Babu and Kesar Singh},
    title     = {On deviations between empirical and quantile
processes for mixing random variables},
    journal   = {Journal of Multivariate Analysis},
    volume    = {8},
    issue     = {4},
    pages     = {532--549},
    year      = {1978},
    month     = {December},
}

@ARTICLE{Berrahou:Bouzebda:Douge:2024,
    author    = {Nour-Eddine Berrahou and Salim Bouzebda and Lahcen Douge},
    title     = {The Bahadur representation for empirical and smooth quantile estimators under association},
    journal   = {Methodology and Computing in Applied Probability},
    volume    = {26},
    number    = {17},
    pages     = {forthgoing},
    year      = {2024},
    month     = {May},
}

@ARTICLE{He:Shao:1996,
    author    = {Xuming He and Qi-Man Shao},
    title     = {A general Bahadur representation of $M$-estimators and its application to linear regression with nonstochastic designs},
    journal   = {The Annals of Statistics},
    volume    = {24},
    number    = {6},
    pages     = {2608--2630},
    year      = {1996},
    month     = {December},
}

@ARTICLE{IMF:2025,
  author       = {{International Monetary Fund}},
  title        = {Rising Financial Stability Risks},
  year         = {2025},
  url          = {https://www.imf.org/external/pubs/ft/ar/2025/in-focus/rising-financial-stability-risks/},
  journal      = {IMF Annual Report 2025 (In Focus)},
  publisher    = {International Monetary Fund}
}

@ARTICLE{Kiefer:1967,
    author    = {Jack Kiefer},
    title     = {On {B}ahadur's representation of sample quantiles},
    journal   = {The Annals of Mathematical Statistics},
    volume    = {38},
    number    = {5},
    pages     = {1323--1342},
    year      = {1967},
    month     = {October},
}

@ARTICLE{Wu:Yu:Wang:2021,
    author    = {Yi Wu and Wei Yu and Xuejun Wang},
    title     = {The Bahadur representation of sample quantiles for $\phi$-mixing random variables and its application},
    journal   = {Statistics},
    volume    = {55},
    number    = {2},
    pages     = {426--444},
    year      = {2021},
    month     = {},
}

\end{document}